\documentclass[11pt]{article}
\PassOptionsToPackage{obeyspaces}{url}
\usepackage[hidelinks]{hyperref}
\usepackage{bookmark}
\bookmarksetup{open,numbered,depth=5}
\usepackage{amsfonts}
\usepackage{amsmath}
\usepackage{amssymb, amsthm, enumitem}
\usepackage{breakurl}
\usepackage{tikz}
\usetikzlibrary{calc}

\usepackage{etoolbox} 
\patchcmd{\thebibliography}{\leftmargin\labelwidth}{\leftmargin\labelwidth\addtolength\itemsep{-0.1\baselineskip}}{}{}

\usepackage{mathtools}
\usepackage{xstring}

\newcommand*\samethanks[1][\value{footnote}]{\footnotemark[#1]}

\author{Boris Bukh\thanks{Department of Mathematical Sciences, Carnegie Mellon University, Pittsburgh, PA 15213, USA\@. Supported in part by U.S.\ taxpayers through NSF CAREER grant DMS-1555149.} \and Zichao Dong\samethanks}

\title{Longest common subsequences between words of very unequal length}
\date{}

\usepackage[nameinlink]{cleveref}

\newtheorem{theorem}{Theorem}
\newtheorem{lemma}[theorem]{Lemma}

\newtheorem{observation}[theorem]{Observation}
\newtheorem{proposition}[theorem]{Proposition}

\theoremstyle{definition}
\newtheorem{game}{Game}

\newcommand*{\eqdef}{\stackrel{\mbox{\normalfont\tiny def}}{=}} 
\newcommand*{\veps}{\varepsilon}                                
\DeclarePairedDelimiter\abs{\lvert}{\rvert}                     

\newcommand*{\wordf}{\textsf}                                   
\newcommand*{\E}{\mathbb{E}}                                    
\DeclareMathOperator{\Var}{Var}                                 
\newcommand*{\N}{\mathbb{N}}                                    
\newcommand*{\Z}{\mathbb{Z}}                                    
\DeclareMathOperator{\LCS}{LCS}                                 
\DeclareMathOperator{\Binom}{Binom}                             
\DeclareMathOperator{\Pois}{Pois}                               
\DeclareMathOperator{\LNDS}{LNDS}                               
\DeclareMathOperator{\Tidy}{T}
\DeclareMathOperator{\Unif}{U}                                  
\DeclareMathOperator{\Penalty}{Pnlt}
\DeclareMathOperator{\len}{len}                                 
\newcommand*{\Strat}{\operatorname{Strat}_0}                    
\DeclareMathOperator{\DStrat}{Strat}                            
\DeclareMathOperator{\Hist}{Hist}                               

\newcommand*{\Med}{m}                                           
\newcommand*{\hl}{d}                                            
\newcommand*{\hlv}{\mathbf{d}}                                  
\newcommand*{\Hl}{D}                                            
\newcommand*{\cA}{\mathcal{A}}                                  
\newcommand*{\starstep}{($\star$) }
\newcommand*{\lmax}{\lambda_{\max}^k}                           
\newcommand*{\wFort}{w_{\text{Fortune}}}                        
\newcommand*{\dTV}{d_{\text{TV}}}                               

\newcommand*{\same}{\text{\textsc{same}}}                       
\newcommand*{\diff}{\text{\textsc{diff}}}                       
\newcommand*{\stin}{\text{\textsc{in}}}
\newcommand*{\stout}{\text{\textsc{out}}}
\newcommand*{\choice}{\text{\textsc{yes}}}
\newcommand*{\nochoice}{\text{\textsc{no}}}

\newcommand*{\Qmin}{Q_{\min}}                                   
\newcommand*{\Qmax}{Q_{\max}}                                   

\newcommand*{\Adv}{I}                                           

\let\oldbegingame\game
\let\oldendgame\endgame
\def\game{\oldbegingame\leavevmode\begin{description}[style=unboxed,nosep,font=\normalfont\itshape,labelindent=1em,leftmargin=!,labelwidth=\widthof{\itshape Game state:}]}
\def\endgame{\end{description}\oldendgame}
\newcommand*{\gamestate}{\item[Game state:]}
\newcommand*{\gamestart}{\item[Start:]}
\newcommand*{\gameturn}{\item[Each turn:]}
\newcommand*{\gameobjective}{\item[Objective:]}
\newcommand*{\gamepenalty}{\item[Penalty:]}
\newcommand*{\gameduration}[1]{\item[Duration:] #1 turns.}

\crefname{enumi}{step}{steps}
\crefname{part}{part}{parts}
\Crefname{game}{Game}{Games}

\begin{document}
\maketitle

\begin{abstract}
We consider the expected length of the longest common subsequence between two random words of lengths $n$ and $(1-\veps)kn$ over a $k$-symbol alphabet.
It is well-known that this quantity is asymptotic to $\gamma_{k,\varepsilon} n$ for some constant $\gamma_{k,\veps}$.
We show that $\gamma_{k,\veps}$ is of the order $1-c\veps^2$ uniformly in $k$ and $\veps$. 
In addition, for large $k$, we give evidence that $\gamma_{k,\veps}$ approaches $1-\tfrac{1}{4}\veps^2$, 
and prove a matching lower bound.
\end{abstract}

\section{Introduction}
\paragraph{Background.} A \emph{word} over alphabet $\Sigma$ is a sequence of elements of $\Sigma$, which we call \emph{symbols}. A \emph{subsequence} in a word $w$ is any word obtained from
$w$ by deleting some, not necessarily contiguous, symbols. By contrast, a \emph{subword} consists of consecutive symbols from $w$. For example,
\wordf{abada} is a subsequence, but not a subword of \wordf{abracadabra}.

For a pair of words $w,w'$, a \emph{common subsequence} is any word that is a subsequence of both $w$ and $w'$. We denote by $\LCS(w,w')$ the length of the 
longest common subsequence between $w$ and $w'$. This quantity is a common way to measure similarity between words. The earliest mathematical studies
of $\LCS$ and its variants were motivated by biological applications \cite{chvatal_sankoff}, but later it found connections to coding theory, linguistics,
and text processing among other fields (book \cite{timewarps} provides general overview, \cite[Ch.~11]{gusfield} discusses computational biology aspects, for recent coding-theoretic applications
see \cite{bukh_guruswami_hastad,synchronization} and references therein).

A particular problem is to understand $\LCS(w,w')$ for a pair of random words $w,w'$. Almost all the work on $\LCS$ for random words concerned
$\LCS(w,w')$ for a pair of random independent \emph{equally long} words $w,w'$. Most of the focus has been on the Chv\'atal--Sankoff constants, which is the limit
\[
  \gamma_k\eqdef \lim \frac{1}{n}\E_{w,w'\sim [k]^n} \LCS(w,w').
\]
Here, we write $\omega\sim \Omega$ to signify that element $\omega$ is chosen uniformly from the set~$\Omega$, with
the convention that whenever this notation occurs several times in the same expression, then the respective choices are independent.
So, $w\sim [k]^n$ signifies that $w$ is a random $n$-symbol word over $[k]$. Similarly, we write $w\sim [k]^{\infty}$ to denote an infinite
random word where each symbol is independently sampled from~$[k]$.

Much work has been done on estimating $\gamma_k$. For $k=2$, Lueker \cite{lueker} proved that $0.788071 \leq \gamma_2\leq 0.826280$. For large $k$,
Kiwi--Loebl--Matou\v{s}ek \cite{kiwi_loebl_matousek} showed that
\begin{equation}\label{eq:kiwi_loebl_matousek}
  \gamma_k=(2+o(1)\bigr)/\sqrt{k}\qquad\text{as }k\to\infty.
\end{equation}
For small $k\geq 3$, the best bounds on $\gamma_k$ can be found in \cite{bgns_chvatal_sankoff_bounds} (upper bounds) and in \cite{kiwi_soto} (lower bounds).
For simulations that estimate $\gamma_k$, see \cite[Table~2]{bgns_chvatal_sankoff_bounds}, \cite[Table~1]{bundschuh2001} and \cite[Section~5]{bukh_cox}.

\paragraph{Results.} In this paper we consider $\LCS$ for words of drastically unequal length. Let us temporarily fix $k\in \mathbb{N}$ and $\veps>0$ arbitrarily. 
It is fairly easy to see that if $w\in [k]^m$ is any fixed word of length $m$, and $w'\sim [k]^{\infty}$, then the shortest prefix of $w'$ that contains
$w$ will be of length about $km$.  So, if $w\sim [k]^n$ and $w'\sim [k]^{(1-\veps)kn}$ are uniform random words of lengths $n$ and $(1-\veps)kn$ respectively,
then $\LCS(w,w')\geq \bigl(1-\veps-o(1)\bigr)n$ with high probability. This bound turns out to be far from being sharp, and our results provide asymptotics for
$\LCS(w,w')$ in this situation.

Let 
\[
  \gamma_{k,\veps}\eqdef \lim_{n\to\infty} \frac{1}{n}\E \LCS(w,w'),\qquad\text{for }w\sim [k]^n,\ w'\sim [k]^{(1-\veps)kn}.
\]
A standard argument using Fekete's lemma and superadditivity of $\LCS$ shows that the limit above exists (see e.g. \cite[Section~1.1]{steele_book}).
\begin{theorem}\label{thm:crude}
For all $k\geq 2$ and all $0<\veps<1/60$,
\[
  1-8\veps^2 \leq \gamma_{k,\veps}\leq 1-\veps^2/72.
\]
\end{theorem}
It is likely that, for fixed $k$ and $\veps\to 0$, the quantity $\gamma_{k,\veps}$ should be asymptotic
to $1-\gamma'_k\veps^2$ for some constant $\gamma_k'$ depending on~$k$. Similarly to the usual Chv\'atal--Sankoff constants, 
determination of constants $\gamma_k'$ for specific values of $k$ appears to be difficult. However, it should be possible
to prove an analogue of \eqref{eq:kiwi_loebl_matousek} for $\gamma_k'$. In fact, we conjecture that
$\gamma_k'\to \tfrac{1}{4}$ as $k\to\infty$. 

We are able to prove a half of this conjecture.
\begin{theorem}\label{thm:asympt}
There exist absolute constants $c,C>0$ such that
$\gamma_{k,\veps}\geq 1-\tfrac{1}{4}\veps^2(1+Ck^{-2/13})$ whenever
$0<\veps\leq c/k \log k$. In particular, if the constants $\gamma_k'$ exist, then $\limsup_{k\to\infty} \gamma_k'\leq \tfrac{1}{4}$.
\end{theorem}
To motivate the conjecture, we must explain the ideas in the proofs of \Cref{thm:crude,thm:asympt}.

\paragraph{Rough proof strategy.}
Both our lower and our upper bounds on $\LCS$ rely on chopping random words $w$ and $w'$ up into linearly many small subwords.
Chopping up is not a new idea, see for example \cite{kiwi_loebl_matousek,bukh_hogenson}, but we do it differently
than previous works. Usually one tries to estimate $\LCS(u,u')$ for a pair of random subwords $u,u'$ of suitably chosen lengths.
In our argument, we chop the longer word $w'$ into subwords of fixed length, but the partition of the shorter word is into subwords
of variable lengths. Namely, for a subword $u'$ of $w'$ and a suffix $u$ of $w$, we shall seek the longest prefix of $u$ that is \emph{almost} a subsequence of $u'$, in the sense that it can
be made into a genuine subsequence by removing only a handful of symbols.

Formally, we say that a word $u$ is \emph{$\hl$\nobreakdash-almost contained} in a word $u'$, and write $u\prec_{\hl} u'$, if we may remove at most $\hl$ symbols from $u$ and obtain a subsequence of $u'$.
In particular, $u\prec_0 u'$ means that $u$ is a subsequence of $u'$.
For example, $\wordf{macabre}\prec_2 \wordf{abracadabra}$. 

We index symbols in a word from $0$, denoting the symbols of a word $u$ by $u[0],u[1],\dotsc$ in order. For example, if $u=\wordf{abad}$, then $u[0]=\wordf{a}$, $u[1]=\wordf{b}$,
$u[2]=\wordf{a}$ and $u[3]=\wordf{d}$. For a word $u$, denote by $u_{<m}$ the prefix of $u$ of length~$m$. 
Consider a pair of random words $w,w'$, where $w\sim [k]^{\infty}$ and $w'$ is a uniform random word of length at least~$L$. Let $P_{\hl}(L)$ be the length
of the longest prefix of $w$ that is $\hl$\nobreakdash-almost contained in~$w'_{<L}$, i.e.,
\[
  P_{\hl}(L)\eqdef \max\{m : w_{<m}\prec_{\hl} w'_{< L}\}.
\]
In other words, if we imagine that $w$ is generated symbol-by-symbol,
$P_{\hl}(L)$ is the waiting time until we obtain a word that is not $\hl$\nobreakdash-almost contained in $w'_{< L}$.

In \Cref{sec:lowerbound,sec:upperbound} we shall derive \Cref{thm:crude} from the following pair of estimates:
\begin{theorem}\label{lem:crudedriftlower}
We have $\E[P_1(L)-P_0(L)]\geq \sqrt{L/7k}$, if $k\geq 2$ and $L\geq 20k$.
\end{theorem}
\begin{theorem}\label{lem:crudedriftupper}
We have $\E[P_{\hl}(L)-P_0(L)]\leq \hl \sqrt{2L/k}+\hl$, for all $k$ and all $\hl$.
\end{theorem}

It is possible to turn tighter estimates on $\E[P_{\hl}(L)-P_0(L)]$ into tighter estimates on $\gamma_{k,\veps}$.
This is precisely how we obtain \Cref{thm:asympt}, by proving the following asymptotics for $\E[P_{\hl}(L)-P_0(L)]$:
\begin{theorem}\label{thm:driftasympt}
We have $\E[P_{\hl}(L)-P_0(L)]=2\sqrt{\hl L/k} \cdot \Bigl(1+O\Bigl(\frac{1}{\hl^{2/3}}+\frac{\log L}{(L/\hl k)^{1/2}}+\frac{\hl^{3/2}}{k^{1/2}}+\frac{k^{1/2}\hl^{3/2}L^{3/2}}{\exp(L^{1/2}k^{-3/2})}\Bigr)\Bigr)$.
\end{theorem}
Sadly, this asymptotics is not sufficiently precise to obtain an upper bound on $\gamma_{k,\veps}$ that
matches the lower bound in \Cref{thm:asympt}. The obstacle is the $\frac{\hl^{3/2}}{k^{1/2}}$ term;
our upper bound arguments require an estimate for $\E[P_{\hl}(L)-P_0(L)]$ be at most linear in~$\hl$. (On the other hand, the last term involving $d^{3/2}$ is unproblematic, as it disappears in the limit $L\to\infty$.)

\paragraph{Connection to the longest non-decreasing subsequences.}
The advantage of focusing on $P_0,\dotsc,P_{\hl}$ is that the growth of $P_{\hl}$ is controlled 
by the length of the longest non-decreasing \mbox{subsequence} (LNDS) in a suitably constructed word over $(\hl+1)$-symbol 
alphabet, for $1\ll \hl\ll k$. 

To explain the connection between the LCS and the LNDS, we must examine
how $P_0(L),\dotsc,P_d(L)$ change as we increase~$L$. 
Recall that $w[i]$ is the symbol of a word $w$ at position $i$, with indexing starting from $0$.
At the start, we have $P_i(0)=i$, for each $i$. We then update $P$'s using the following observation.
\begin{proposition}[Proof is in \Cref{sec:props}]\label{prop:update}
We may compute the values of $P_0(L+1),\dotsc,P_d(L+1)$ from $P_0(L),\dotsc,P_d(L)$ by doing the following, in order:
\begin{enumerate}[label=(\Alph*), ref=(\Alph*), series=updateseries]
\item \label{step:Achoice} Examine $w'[L]$. Set $A[L]=\{ i : w[P_i(L)]=w'[L] \}$.
\item \label{step:Aadvance}For each $i=0,1,\dotsc,d$ in order, do
\begin{itemize}
\item If $i\in A[L]$, then $P_i(L+1)=P_i(L)+1$.
\item If $i\notin A[L]$, then $P_i(L+1)=\max(P_i(L),P_{i-1}(L+1)+1)$. (With the convention that $P_{-1}(L)=-\infty$.)
\end{itemize}
\end{enumerate}
\end{proposition}
We may describe this process alternatively by imagining that $P_0,\dotsc,P_{\hl}$ indicate 
positions of $\hl+1$ particles. When we expose the value of $w'[L]$, we 
check which particles sit atop matching symbols; we denote those by $A[L]$.
We then advance particles in $A[L]$ one step to the right. If two particles
collide, the particle on the left `bumps' the particle on the right, causing it
to advance one step further to the right.

Below is an example of the first three exposure steps for words $w=\wordf{1323121\ldots}$ and $w'=\wordf{231\ldots}$.
In this example, we write the elements of sets $A[0],\dotsc,A[L]$ in descending order.

\begin{center}
\def\symdist{0.33} 
\def\stepdist{4.4}
\newcommand*{\basepic}[8]{
\pgfmathsetmacro{\xshft}{#1*\stepdist}
\begin{scope}[xshift=\xshft cm]
\node at (0.5,1) {$\scriptstyle w$};
\foreach \i / \s in {0/1,1/3,2/2,3/3,4/1,5/2,6/1,7/\ldots}  
  \node at ($(1.0,1)+\i*(\symdist,0)$) {$\scriptstyle \wordf{\s\vphantom{?}}$};
\node at (0.5,0) {$\scriptstyle w'$};
\foreach \i / \s in {0/#6,1/#7,2/#8,3/?,4/?,5/?,6/?,7/\ldots}            
  \node at ($(1.0,0)+\i*(\symdist,0)$) {$\scriptstyle \wordf{\s\vphantom{?}}$};
\foreach \i / \p in {#2/0,#3/1,#4/2,#5/3} 
{
  \pgfmathsetmacro{\x}{1.0+\i*\symdist}
  \node at (\x,1.8) {$\scriptstyle P_{\p}$};
  \draw[->] (\x,1.6) -- (\x,1.2);
}
\draw[very thin, rounded corners=2pt] (0.25,-0.25) rectangle ($(1.2,2)+7*(\symdist,0)$);
\node at ($(0.475,-0.5)+3.5*(\symdist,0)$) {$\scriptstyle L=#1$};
\end{scope}
}
\newcommand*{\transpic}[2]{%
\pgfmathsetmacro{\x}{0.225+3.5*\symdist+(#1+0.5)*\stepdist}
\draw[->] (\x,0.65) -- +(1,0); 
\node at ($(\x,0.85)+(0.45,0)$) {$\scriptscriptstyle A[#1]=\wordf{#2}$};
}
\begin{tikzpicture}
\basepic{0}{0}{1}{2}{3}{?}{?}{?}
\transpic{0}{2}
\basepic{1}{0}{1}{3}{4}{2}{?}{?}
\transpic{1}{21}
\basepic{2}{0}{2}{4}{5}{2}{3}{?}
\transpic{2}{20}
\basepic{3}{1}{2}{5}{6}{2}{3}{1}
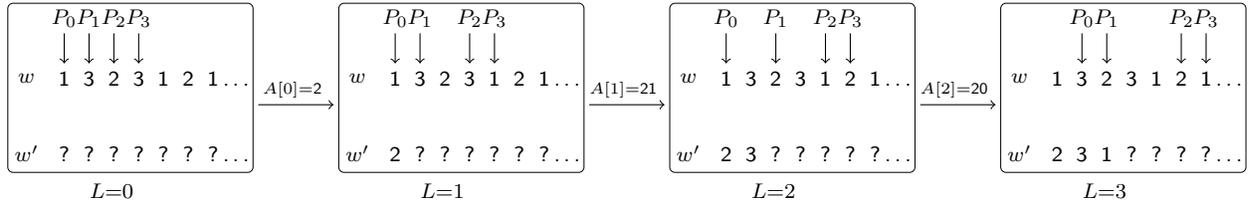
\node at ($(0.475,-1.15)+3.5*(\symdist,0)+1.5*(\stepdist,0)$) {\textbf{\hypertarget{figone}{Figure 1}: }Evolution of $P_0,P_1,P_2,P_3$ (example).};
\end{tikzpicture}
\end{center}

The values of $P_0(L),\dotsc,P_d(L)$ depend only on sets $A[0],\dotsc,A[L-1]$. We can describe this dependence in
terms of non-decreasing subsequences. A word $w$ is \emph{non-decreasing} if $w[i]\leq w[i+1]$ holds for all $i$. 
For example, the word \wordf{001224} is non-decreasing.
For a word $w$, let $\LNDS(w)$ denote the length of the longest non-decreasing subsequence of $w$.
Slightly more generally, let $\LNDS_i(w)=\nobreak\LNDS(w|_i)$, where $w|_i$ is the word
obtained by deleting symbols greater than $i$ from $w$. In particular,
if $w\in \{0,1,\dotsc,d\}^L$, then $\LNDS_d(w)=\LNDS(w)$.

For a set $A\subseteq \{0,1,\dotsc,d\}$, consider the word obtained by writing the elements of~$A$
in decreasing order. We denote this word by the same letter~$A$. 
\begin{proposition}[Proof is in \Cref{sec:props}]\label{prop:lcstolnds}
Let $A[0],A[1],\dotsc,A[L-1]$ be sets as in \Cref{prop:update}. Then
the word $A[0]A[1]\dotsb A[L-1]$ satisfies $\LNDS_i(A[0]A[1]\dotsb A[L-1])=P_i(L)-i$, for every~$i$.
\end{proposition}
For example in the \hyperlink{figone}{Figure~1}, we have $P_3(3)=6$ and $\LNDS(\wordf{22120})=3$. We
also have $P_1(3)=2$ and $\LNDS_1(\wordf{22120})=1$ in the same example.

\paragraph{Expectant partitions.} Consider the update rule in \Cref{prop:update}.
Let $\cA[L]$ consist of all non-empty sets of the form $\{i: w[P_i(L)]=s\}$ for $s\in [k]$.
Note that $\cA[L]$ is a partition of $\{0,1,\dotsc,d\}$.
If the word $w'$ is random, then the set $A[L]$ at \Cref{step:Achoice} in \Cref{prop:update}
is chosen uniformly from $\cA[L]$, conditioned on $A[L]$ being non-empty.
For that reason, we call $\cA[L]$ \emph{expectant partition} at step~$L$. For example,
in \hyperlink{figone}{Figure~1} the expectant partitions are 
$\cA[0]=\bigl\{\{0\},\{1,3\},\{2\}\bigr\}$, $\cA[1]=\bigl\{\{0,3\},\{1,2\}\bigr\}$, $\cA[2]=\bigl\{\{0,2\},\{1,3\}\bigr\}$, and $\cA[3]=\bigl\{\{0\},\{1,2\},\{3\}\bigr\}$.
At the next step, the set $A[3]$ will be chosen to be one of $\{0\}$, $\{1,2\}$, $\{3\}$, each with probability $\tfrac{1}{3}$.

We call the partition all of whose sets are singletons \emph{trivial partition}.
Had all expectant partitions been trivial, then $A[0]A[1]\dotsb A[L-1]$ would have been
a uniform random word. The behavior of LNDS on uniform random words is well understood, thanks to the
work of Tracy and Widom \cite{tracy_widom} (see also \cite{johansson,kuperberg}).

In our proof of \Cref{thm:asympt}, we will first show that most expectant partitions are trivial.
Then we will show that the remaining handful of non-trivial partitions do not change $P_d(L)$ much. 

\paragraph{Adversarial game arguments.} Our main technical innovation concerns 
analysis of Markov chains. The Markov chains that arise in our analysis
of $P_d(L)$ have complicated state spaces and complicated transition rules.
Instead of trying to describe their behavior directly, we choose to 
ignore certain details of the chain, and do the worst-case analysis instead.

More formally, we imagine that certain transitions are no longer random, but instead
are chosen by a suitably restricted adversary. Every adversary's strategy leads to a different
Markov chain, with one of the choices being our original chain. Furthermore,
we can reduce the state space of the chain, by ignoring those parts that are
under adversary's control. When this is done carefully,
we are able to ensure that the adversary's optimal strategy 
leads to a much smaller Markov chain that is easier-to-analyze 
than the original chain.

We shall formalize these ideas in \Cref{lem:adversarial}.

\paragraph{Paper organization.}
We begin by showing how to turn the estimates on $\E[P_d-P_0]$ into bounds on $\gamma_{k,\veps}$:
In \Cref{sec:lowerbound} we derive the lower bounds in \Cref{thm:crude,thm:asympt} from the lower bounds on $\E[P_d-P_0]$
in \Cref{lem:crudedriftlower,thm:driftasympt}. Similarly, in \Cref{sec:upperbound} we derive the upper
bound in \Cref{thm:crude} from the upper bound on $\E[P_d-P_0]$ in \Cref{lem:crudedriftupper}.

The main bulk of paper is then devoted to proving bounds on $\E[P_d-P_0]$. 
Since these bounds all use the adversarial game argument, we start by formalizing the 
argument in \Cref{sec:adversarial}. We then present the proofs of \Cref{lem:crudedriftlower,lem:crudedriftupper},
starting with the easier \Cref{lem:crudedriftupper} in \Cref{sec:crudedriftupper},
and following it up with the proof of \Cref{lem:crudedriftlower} in \Cref{sec:crudedriftlower}.

The asymptotic result in \Cref{thm:driftasympt} is broken into several parts. 
In \Cref{sec:trivial} we show that most expectant partitions are trivial. We 
then estimate the effect of non-trivial partitions in \Cref{sec:asympt}. 

\paragraph{Acknowledgments.} We are grateful to Greg Kuperberg and Kurt Johansson for answering our queries
about longest non-decreasing subsequences. We also benefited from discussions with Chris Cox and Zilin Jiang.
We are very thankful to the referee for extensive comments on the manuscript.

\section{Proof of the lower bounds in \texorpdfstring{\Cref{thm:crude,thm:asympt}}{Theorems~\ref{thm:crude} and~\ref{thm:asympt}}}\label{sec:lowerbound}
In this section, we will show the following:
\begin{lemma}\label{lem:lowerbnd}
Suppose $\hl,k\in \N$ and $\alpha>0$ are such that $\E[P_{\hl}(L)-P_0(L)]\geq \alpha \sqrt{L/k}$ holds for all $L\geq L_0$.
Then $\gamma_{k,\veps}\geq 1-d\veps^2(1+4\veps)/\alpha^2$ for all $0<\veps<\min(1/20,\alpha \sqrt{k/2L_0})$.
\end{lemma}
Upon applying the lemma with $d=1$, $\alpha=1/\sqrt{7}$ and $L_0=200k$, the lower bound in \Cref{thm:crude} instantly follows from \Cref{lem:crudedriftlower}.
Similarly, to derive \Cref{thm:asympt} from \Cref{thm:driftasympt} we apply the lemma with $d=k^{3/13}$ and
$L_0=36 k^3\log^2 k$.\medskip

With hindsight let 
\[
  L\eqdef (1-2\veps)^2\alpha^2 k \veps^{-2}, \qquad M\eqdef (1-\veps)kn/L.
\]

We employ the customary abuse of notation: we write our proofs as if the numbers $L$
and $M$ are integers. This can be made formal by rounding $L$ to an integer,
and truncating the word of length $(1-\veps)kn$ slightly so that its
length is an integral multiple of $L$. Doing so does not affect the limit $\gamma_{k,\veps}$,
which is the subject of the present lemma. We will use similar abuses of notation later in the
paper without any further comment.\smallskip

Let $w\sim [k]^n$ and $w'\sim [k]^{(1-\veps)kn}$ be random words as in the definition of $\gamma_{k,\veps}$.
Imagine that the word $w$ is a prefix of an infinite random word $\overline{w}\sim [k]^{\infty}$.
Write $w'$ as a concatenation of $M$ words of length $L$ each, say $w'=w_1'w_2'\dotsb w_M'$.
We iteratively define words $\overline{w}_1,\dotsc,\overline{w}_M$ in such a way that
$\overline{w}_1\dotsb\overline{w}_M$ is a prefix of $\overline{w}$: 
Given $\overline{w}_1$ through $\overline{w}_r$, write $\overline{w}=\overline{w}_1\dotsb\overline{w}_r v$ for some suffix~$v$.
We then define $\overline{w}_{r+1}$ to be the longest prefix of $v$ that is $\hl$\nobreakdash-almost contained in~$w'_{r+1}$.

Let $Y_i\eqdef \len \overline{w}_i$. The random variables $Y_1,\dotsc,Y_M$ are independent. Indeed,
we may imagine first generating $w_1'$ and $\overline{w}_1$, then generating $w_2'$ and $\overline{w}_2$, and so forth.
In this process, to generate the pair $(w_{\ell}',\overline{w}_{\ell})$, we first generate $w_{\ell}'$, and 
then we append symbols to $\overline{w}_{\ell}$ as long as $\overline{w}_{\ell} \prec_d w_{\ell}'$.
Hence, after the pair $(w_{\ell}',\overline{w}_{\ell})$ has been generated, the symbol following
$\overline{w}_{\ell}$ has also been generated. This symbol will become the first symbol
of $\overline{w}_{\ell+1}$, and so the word $\overline{w}_{\ell+1}$ will not be independent
from the words $\overline{w}_1,\dotsc,\overline{w}_{\ell}$. However,
since the first symbol of $w_{\ell+1}'$ is independent from $w_1',\dotsc,w_{\ell}',\overline{w}_1,\dotsc,\overline{w}_{\ell},\overline{w}_{\ell+1}[0]$,
the event $w_{\ell+1}'[0]=\overline{w}_{\ell+1}[0]$ still has probability $1/k$, and so
the \emph{length} of $\overline{w}_{\ell+1}$ is independent from the lengths of the preceding words.

Let $Y\eqdef Y_1+Y_2+\dotsb+Y_M$.
If $Y\geq n$, then the word $w$ is a prefix of $\overline{w}_1\dotsb\overline{w}_M$,
implying that $\LCS(w,w')\geq n-dM$.

Since the lengths of $\overline{w}_1,\dotsc,\overline{w}_M$ are independent, 
$Y=Y_1+\dotsb+Y_M$ is a sum of independent random variables. Because each $Y_i$ is sampled from~$P_d(L)$, and
$\E[P_0(L)]=L/k$, it follows that
\[
  \E[Y]\geq M(L/k+\alpha\sqrt{L/k})=\frac{(1-\veps)^2}{1-2\veps}n\geq (1+\veps^2)n.
\]
Since $d\leq Y_i\leq L+d$ for each $i$, the Chernoff bound \cite[Theorem~A.1.18]{alon_spencer} implies that
\[
  \Pr[Y\leq n]\leq \Pr[Y\leq \E Y-\veps^2 n]\leq \exp\bigl(-2(\veps^2 n)^2/ML^2\bigr)\leq \exp\bigl(-2\veps^4 n/kL\bigr)=o(1).
\]
Hence,
\[
  \E \LCS(w,w')\geq \Pr[Y\geq n](n-dM)\geq \Bigl(1-d\frac{1-\veps}{(1-2\veps)^2}\cdot \frac{\veps^2}{\alpha^2}-o(1)\Bigr)n.
\]
As $\frac{1-\veps}{(1-2\veps)^2}\leq 1+4\veps$ for $\veps\leq 1/20$, the proof is complete.

\section{Proof of the upper bound in \texorpdfstring{\Cref{thm:crude}}{Theorem \ref{thm:crude}} assuming \texorpdfstring{\Cref{lem:crudedriftupper}}{Theorem \ref{lem:crudedriftupper}}}
\label{sec:upperbound}
In the introduction, we defined a subsequence of $w$ as any word obtained from $w$ by removing some of the symbols.
Sometimes the same word can be so obtained from $w$ in several ways. To eliminate this ambiguity, we introduce a couple
of definitions.

Recall that $w[i]$ denotes the symbol of a word $w$ at position $i$, with indexing starting from $0$.
For a word $w$ and a set of integers $I$, we define $w[I]\eqdef (w[i] : i\in I)$. For example, if $w=\wordf{abracadabra}$
and $I=\{1,2,4\}$, then $w[I]=\wordf{brc}$. A \emph{distinguished subsequence} of a word $w\in [k]^n$ is a pair $(u,I)$ 
such that $u=w[I]$; note that $u$ is the usual (undistinguished) subsequence of~$w$. 
Similarly, a \emph{distinguished common subsequence}
of words $w$ and $w'$ is a triple $(u,I,I')$ such that $u=w[I]=w'[I']$. When there is no chance of confusion, 
we shall use $u$ to refer to the distinguished (common) subsequence, omitting $I$ and $I'$ from the notation.

For each pair of words $w\in [k]^n$ and $w'\in [k]^{(1-\veps)kn}$, we shall select a suitable distinguished longest common subsequence, and assign a vector $\hlv$
that describes the ``shape'' of that common subsequence. We will then show that, for any fixed $\hlv$, a pair of random words is very
unlikely to have a long common subsequence of that shape. The union bound will then complete the argument.

\paragraph{Standard prefix.} Define three constants 
\[
  L\eqdef 432 k\veps^{-2},\qquad M\eqdef (1-\veps)kn/L,\qquad \Hl\eqdef \tfrac{1}{2}\veps n\sqrt{k/3L}.
\]
For the future use, observe that the choice of $L$ implies that $\sqrt{2L/k}+1\leq \sqrt{3L/k}$.\smallskip

As in the preceding section, whenever we have a word $w'\in [k]^{(1-\veps)kn}$ of length $(1-\veps)kn$, we write it 
as a concatenation of $M$ words of length $L$ each, $w'=\nobreak w_1'\dotsb w_M'$. 
We use this notation throughout the section.

Suppose $\hlv=(\hl_1,\dotsc,\hl_M)$ is a vector of non-negative integers.
Given $w\in [k]^n$ and $w'\in [k]^{(1-\veps)kn}$, we say that $w$ is \emph{$\hlv$\nobreakdash-almost contained} in $w'$
and write $w\prec_{\hlv} w'$ if there is a decomposition of $w$ as $w=w_1\dotsb w_M$
for some words $w_1,\dotsc,w_M$ such that the word $w_i$ is a $\hl_i$\nobreakdash-almost contained in~$w_i'$.

\begin{observation}\label{obs:domination}
If $\LCS(w,w')\geq n-\Hl$, then there is a vector $\hlv=(\hl_1,\dotsc,\hl_M)\in \Z_+^M$ satisfying
$\hl_1+\dotsb+\hl_M\leq \Hl$ such that $w\prec_{\hlv} w'$.
\end{observation}
\begin{proof}
Let $u$ be a distinguished common subsequence between $w$ and $w'$ of length $n-\Hl$.
Define $w_1,\dotsc,w_M$ inductively. Assume that $w_1,\dotsc,w_{i-1}$ have been defined, and
$w=w_1\dotsb w_{i-1}v$ for some word $v$. Let $w_i$ be the longest prefix of $v$
that contains no symbol that $u$ matches to $w_{i+1}'$. Let $\hl_i$ be the number of unmatched symbols in $w_i$.
It is then clear that $\len u+\sum \hl_i= n$.
\end{proof}

Let a vector $\hlv=(\hl_1,\dotsc,\hl_M)$ and words $w\in [k]^n$ and $w'\in[k]^{(1-\veps)kn}$ be given. Define a sequence of words $w_1,\dotsc,w_M$ inductively as follows:
Assume that $w_1$ through $w_{i-1}$ have been defined, and write $w=w_1\dotsb w_{i-1}v$. Let $w_i$ be the longest prefix of $v$ that is $\hl_i$\nobreakdash-almost contained in $w_i'$. It might happen that,
for some $j$, the word $v$ is $\hl_j$\nobreakdash-almost contained in $w_j'$. In that case, $w_j=v$
and so all the subsequent words $w_{j+1},w_{j+2},\dotsc,w_M$ are empty.
We call $w_1\dotsb w_M$ the \emph{standard prefix} of $w$ (with respect to $w'$ and~$\hlv$). 

\begin{lemma}\label{lem:stdprefix}
If $w\prec_{\hlv} w'$, then the standard prefix of $w$ with respect to $w'$ and $\hlv$ is equal to $w$.
\end{lemma}
\begin{proof}
Let $w=\hat{w}_1\dotsb \hat{w}_M$ be a decomposition of $w$ satisfying
$\hat{w}_i\prec_{\hl_i} w_i'$, which exists because $w\prec_{\hlv} w'$.
Let $w_1\dotsb w_M$ be the standard prefix of $w$. By induction on $i$ it follows that
$\hat{w}_1\dotsb \hat{w}_i$ is a prefix of $w_1\dotsb w_i$. In particular,
$w=\hat{w}_1\dotsb \hat{w}_M$ is a prefix of $w_1\dotsb w_M$. So, $w=w_1\dotsb w_M$ follows.
\end{proof}

\paragraph{Any fixed \texorpdfstring{$\hlv$-vector}{d-vector} is unlikely.} Let $\hlv=(\hl_1,\dotsc,\hl_M)$ be an arbitrary vector
of nonnegative integers satisfying $\hl_1+\dotsb+\hl_M\leq \Hl$. We shall estimate
the probability that $w\prec_{\hlv} w'$ for random $w\sim [k]^n$ and $w'\sim [k]^{(1-\veps)kn}$.\medskip

We imagine generating symbols of $w$ and $w'$ gradually in $M+1$ \emph{rounds}, numbered $0,1,\dotsc,M$. 
In round $0$ we generate only the first symbol of~$w$. Before the start of round $i$ we already
generated 
\begin{itemize}
\item subwords $w_1',\dotsc,w_{i-1}'$ of $w'$, and
\item the first $i-1$ words $w_1,\dotsc,w_{i-1}$ in the standard prefix of $w$, and
\item a single symbol that follows $w_{i-1}$ in $w$.
\end{itemize}
In round $i$, we first generate $L$ symbols that make up $w_i'$. We then
generate symbols of $w_i$ (starting with a single already-generated symbol)
as long as $w_i\prec_{\hl_i} w_i'$. 
Once $w_i \prec_{\hl_i} w_i'$ ceases to hold, we backtrack one symbol. That symbol will become the first symbol of $w_{i+1}$ in the next round. 
(This backtracking is similar to the backtracking in \Cref{sec:lowerbound} in the definition of $Y_i$.)

We also need a pair of \emph{fix-up rules}: once we generated $n$ symbols of $w$, 
we stop generating symbols of $w$, and generate all the remaining symbols of $w'$ in one shot. 
Similarly, if after the end of the $M$'th round, the word $w$ still has fewer than $n$ symbols, 
we generate the remaining symbols of~$w$. \medskip

It is clear that this algorithm generates a pair of independent uniformly
distributed words $w\sim [k]^n$ and $w'\sim [k]^{(1-\veps)kn}$. To simplify
the analysis, it is convenient to consider the version of this algorithm
without the fix-up rules. We call this modified algorithm \emph{tidy},
and denote the distribution on pairs of words $(w,w')$ that it induces by $\Tidy(\hlv)$.

Let $\Unif$ denote the uniform distribution on $[k]^n\times [k]^{(1-\veps)kn}$.
Since, by \Cref{lem:stdprefix}, $w\prec_{\hlv} w'$ implies that $\len w_1+\dotsb+\len w_M=\len w$,
it follows that
\begin{equation}\label{eq:tidycoupling}
  \Pr_{(w,w')\sim \Unif}[w\prec_{\hlv} w']=\Pr_{(w,w')\sim \Tidy(\hlv)}[\len w \geq n]=\Pr_{(w,w')\sim \Tidy(\hlv)}[\len w_1+\dotsb+\len w_M\geq n].
\end{equation}

Because the first symbol of $w_i$ depends on the first $i-1$ rounds, 
the word $w_i$ is \emph{not} independent of $w_1$ through $w_{i-1}$ and $w_1'$ through $w_{i-1}'$.
However, the distribution of $\len w_i$ is the same for every fixed initial symbol.
So, the lengths $\len w_1,\dotsc,\len w_M$ are independent.

We shall use Talagrand's inequality to bound the probability on the right side of \eqref{eq:tidycoupling}.
Recall that, in the context of Talagrand's inequality, a random variable $W$ on a product space $\Omega$ is called \emph{$f$\nobreakdash-certifiable}
if, whenever, $W(x)\geq b$, there exists a set of at most $f(b)$ coordinates such that every $y\in\Omega$ agreeing with
$x$ on these coordinates satisfies $W(y)\geq b$.

Sample $w'\sim [k]^{(1-\veps)kn}$ and $M$ independent infinite words $w^{(1)},\dotsc,w^{(M)}\sim [k]^{\infty}$.
Let $\hat{w}_i$ be the longest prefix of $w^{(i)}$ satisfying $\hat{w}_i\prec_{\hl_i} w_i'$. The vectors
$(\len w_1,\dotsb,\len w_M)$ and $(\len \hat{w}_1,\dotsb,\len \hat{w}_M)$ are identically distributed.
Define a random variable $Y\eqdef \len \hat{w}_1+\dotsb+\len \hat{w}_M$. 
Since $\len \hat{w}_i$ is sampled from the distribution $P_{\hl_i}(L)$, from \Cref{lem:crudedriftupper}
it follows that
\begin{align}
 \E[Y]&\geq M\E[P_0(L)]=ML/k=(1-\veps)n,\label{eq:lowerew}\\
 \E[Y]&\leq \sum_{i=1}^M \bigl(L/k+ \hl_i\sqrt{2L/k}+\hl_i\bigr)\leq ML/k+\Hl(\sqrt{2L/k}+1)\notag\\&\leq ML/k+\Hl\sqrt{3L/k}=(1-\veps/2)n.\label{eq:upperew}
\end{align}
The random variable $Y= \len \hat{w}_1+\dotsb+\len \hat{w}_M$ is $2b$\nobreakdash-certifiable. Indeed, let $\Omega$ be
the space of all $(M+1)$-tuples $(w^{(1)},\dotsc,w^{(M)},w')\in ([k]^{\infty})^M\times [k]^{(1-\veps)kn}$ endowed with the uniform measure.
We may think of $\Omega$ as a product of infinitely many uniform random variables sampled from~$[k]$ (one for each symbol of $w^{(1)},\dotsc,w^{(M)},w'$). 
If $Y(w^{(1)},\dotsc,w^{(M)},w')\geq b$, then there are $b$ symbols in $w'$ and a total of $b$ symbols in words $w^{(1)},\dotsc,w^{(M)},w'$ that make up $\hat{w}_1$ through~$\hat{w}_M$.
These $2b$ symbols certify that $Y\geq b$. Hence, Talagrand's inequality 
(see \cite[Theorem~7.7.1]{alon_spencer}) tells us that
\begin{equation}
  \Pr\bigl[Y\leq b-t\sqrt{2b}\bigr]\Pr[Y\geq b]\leq \exp(-t^2/4).\label{eq:talagrand}
\end{equation}
Let $\Med$ be the median of~$Y$. We apply \eqref{eq:talagrand} with $b=\Med$
and with $b=\Med+t^2+2t\sqrt{2\Med}$. The first choice gives us 
\begin{align}
  \Pr[Y\leq \Med-t\sqrt{2\Med}]&\leq 2\exp(-t^2/4).\label{eq:lowertalagrand}\\
\intertext{Because of $\Med+t^2+t\sqrt{2\Med} -t\sqrt{2(\Med+t^2+t\sqrt{2\Med})}\geq \Med$, the second choice
gives us}
\Pr[Y\geq \Med+t^2+t\sqrt{2\Med}]&\leq 2\exp(-t^2/4).\label{eq:uppertalagrand}
\end{align}
From \eqref{eq:lowertalagrand} we deduce that
\begin{align*}
\E[Y]\geq \Med-\sqrt{2\Med}\sum_{t=0}^{\infty} \Pr[Y\leq \Med-t\sqrt{2\Med}]\geq \Med-2\sqrt{2\Med}\sum_{t=0}^{\infty}\exp(-t^2/4)\geq \Med-5\sqrt{2\Med}.
\end{align*}
Note that this implies that $\Med\leq 2n$, for otherwise $\E[Y]>n$ contradicting \eqref{eq:upperew}. So, \eqref{eq:uppertalagrand} implies


\begin{equation}\label{eq:talagrandfinal}
\begin{aligned}
 \Pr[Y\geq \E[Y]+t^2+2(t+5)\sqrt{n}]
&\leq  \Pr[Y\geq \E[Y]+t^2+(t+5)\sqrt{2\Med}]\\
&\leq \Pr[Y\geq \Med+t^2+t\sqrt{2\Med}]\\&\leq 2\exp(-t^2/4).
\end{aligned}
\end{equation}
We choose $t=\tfrac{\veps}{6}\sqrt{n}$. 
With this choice, $t^2+2(t+5)\sqrt{n}<3t\sqrt{n}\leq \veps n/2$, and so
from the combination of \eqref{eq:tidycoupling}, \eqref{eq:upperew}, and \eqref{eq:talagrandfinal} we obtain
\begin{equation}\label{eq:unlikely}
 \Pr_{(w,w')\sim \Unif}[w\prec_{\hlv} w']\leq 2\exp(-\veps^2 n/144).
\end{equation}


\paragraph{Putting everything together.} 
The number of nonnegative integer vectors $(\hl_1,\dotsc,\hl_M)$ satisfying $\hl_1+\dotsb+\hl_M\leq \Hl$
is $\binom{\Hl+M}{M}$. Plugging in our choice of constants $\Hl$ and $M$, and using the bound $\binom{m}{cm}\leq 2^{H(c)m}$, 
we obtain
\[
  \binom{\Hl+M}{M}=\binom{\veps^2 n(1/72+1/432)}{\veps^2 n/432}\leq \exp(\veps^2 n/150).
\]
\Cref{obs:domination} and inequality \eqref{eq:unlikely} then tell us that
\[
  \Pr[\LCS(w,w')\geq n-\Hl]\leq \binom{\Hl+M}{M}\cdot 2\exp(-\veps^2 n/144)=o(1/n).
\]
As this implies that $\E \LCS(w,w')\leq n-\Hl-1$, the proof is complete.

\section{Adversarial game argument}\label{sec:adversarial}
Suppose $X_0,X_1,X_2,\dotsc$ is a Markov chain with a complicated transition rule, and that our goal is to 
estimate $\E f(X_n)$. One approach is to replace some of the randomness in the chain
by adversarial choices. We can then bound $\E f(X_n)$ from below by $\inf \E f(X_n)$, where
the infimum is over all possible strategies for the adversary. 

Furthermore, we may do so in stages: we first replace some of the randomness by adversarial choices,
determine what optimal (or near-optimal) choices are, use that to simplify the chain, and then replace
some of the remaining randomness by adversarial choices again.

The aim of this section is to make these ideas precise. 
The reader might want to first glance at the applications in the later sections, and
only then continue reading the formalism that underlies them.

\paragraph{Markov games.}
To make these ideas precise, we need to define the notion of an `adversarial strategy'.
Formally, we define a \emph{Markov game} as a tuple $(\Adv,\Omega,P,R)$, where
\begin{itemize}
\item The set $\Adv\subseteq \{0,1,\dotsc\}$ indexes the times in which the transitions are chosen by the adversary,
\item $\Omega_0,\Omega_1,\Omega_2,\dotsc$ is a sequence of at most countable sets, with $\abs{\Omega_0}=1$, where $\Omega_n$ represents the possible states at the $n$'th step,
\item $R$ is a tuple of functions $R=(R_n)_{n\in \Adv}$ with $R_n\colon \Omega_n\to 2^{\Omega_{n+1}}\setminus\{\emptyset\}$, where $R_n(x_n)\subseteq \Omega_{n+1}$ represents the set of states to which the adversary can move from the state $x_n\in \Omega_n$,
\item $P$ is a tuple of functions $P=(P_n)_{n\notin \Adv}$ with $P_n\colon \Omega_{n+1}\times \Omega_n\to [0,1]$ that satisfies \linebreak$\sum_{x_{n+1}\in \Omega_{n+1}} P_n(x_{n+1},x_n)=1$ for each $x_n\in\Omega_n$. The value of $P_n(x_{n+1},x_n)$ represents 
the transition probability to the state $x_{n+1}\in \Omega_{n+1}$ from the state $x_n\in \Omega_n$.
\end{itemize}
In this language, \emph{Markov chain}\footnote{In this paper, Markov chains are not assumed to be time-homogeneous, i.e., we allow the state spaces and the transition probabilities to depend
on the time step $n$.} is the same as a Markov game with $\Adv=\emptyset$. The only element in $\Omega_0$ is the initial state of the Markov chain.
Since the tuple $R$ is trivial if $I=\emptyset$, we abbreviate $(\emptyset,\Omega,P,R)$ to $(\Omega,P)$ in a case of a Markov chain.

Note that a Markov game is \emph{not} a probability space, and is not naturally associated with a probability space (unless $I=\emptyset$). Instead it models
a one-player game in which some of the turns are made by the player (which we call Adversary), and some turns are made at random.
When discussing the applications of actual Markov games in the following sections, for ease of reading we shall usually personify the randomness as Fortune.
In this section, however, our task is to lay the formal basis for the adversarial game argument, and so we use the formal definition directly.\medskip

The next definition captures the notion of `giving the adversary more choices'. 
Suppose that $G=\nobreak(\Adv,\Omega,P,R)$ and $G'=(\Adv',\Omega',P',R')$ are a pair of Markov games. We say that a sequence $\pi=\nobreak(\pi_0,\pi_1,\dotsc)$
is a \emph{morphism} from $G$ to $G'$ if
\begin{enumerate}[label=(M\arabic*), ref=(M\arabic*)]
\item $\Adv\subseteq \Adv'$,
\item $\pi_n$ is a function from $\Omega_n$ to $\Omega_n'$, for each $n$,
\item \label{morphism:inclusion} for each $n\in \Adv$, and $x_n\in \Omega_n$, we have $\pi_{n+1}(R_n(x_n))\subseteq R_n'(\pi_n(x_n))$,
\item \label{morphism:mild} for each $n\in \Adv'\setminus \Adv$, whenever $x_n\in \Omega_n$ and $x_{n+1}\in \Omega_{n+1}$ satisfy
$P_n(x_{n+1},x_n)>0$, the set $R_n'(\pi_n(x_n))$ contains $\pi_{n+1}(x_{n+1})$,
\item \label{morphism:strong} for each $n\notin \Adv'$, every $x_n\in \Omega_n$ and $x_{n+1}'\in \Omega_{n+1}'$ satisfy
\[
  P_n'\bigl(x_{n+1}',\pi_n(x_n)\bigr)=\sum_{x: \pi_{n+1}(x)=x_{n+1}' } P_n(x,x_n).
\]
\end{enumerate}
If $\Omega_n=\Omega_n'$ and each $\pi_n$ is the identity function, then this definition perfectly mirrors the intuitive
idea of giving more choices to the adversary. In general, if $\pi_n$ is not injective, say $\pi_n(x_n)=\pi_n(\tilde{x}_n)$, then 
the morphism $\pi$ collapses the states $x_n$ and $\tilde{x}_n$ into a single state. The conditions
\ref{morphism:mild} and \ref{morphism:strong} require that, in the Markov game $G'$, the new state is indistinguishable from its preimages in~$G$.

\paragraph{Strategies.} Intuitively, a deterministic strategy for the adversary
consists of choosing, for every time step $n\in \Adv$ and every state $x_n\in \Omega_n$, 
a state in $R_n(x_n)$ that $x_n$ is mapped to. A probabilistic strategy is then identified
with a collection of probability distributions on all $R_n(x_n)$, which is quite an unwieldy
object to work with. 

We adopt a simpler point of view. Any choice of a (possibly probabilistic) strategy 
turns a Markov game into a Markov chain. 
In our approach, we identify the strategy with this Markov chain. 
Formally,
we say that a Markov chain $S=(\Omega,P)$ is a \emph{strategy} for $G'=(\Adv',\Omega',P',R')$ if there is a morphism $\pi\colon S\to G'$ such that
\begin{itemize}
\item $\Omega_n=\Omega_n'$, and
\item there is a morphism from $S$ to $G'$ in which each $\pi_n\colon \Omega_n\to\Omega_n'$ is the identity map.
\end{itemize}
Let $\Strat(G')$ denote the set of all strategies for $G'$.

Say that a Markov chain $S=(\Omega,P)$ is \emph{determined at $n$} if $P_n(x_{n+1},x_n)\in \{0,1\}$ for all $x_n\in \Omega_n$.
We say that a strategy for a Markov game $G=(\Adv,\Omega,P,R)$ is \emph{deterministic} if it is determined for every $n\in \Adv$.
Let $\DStrat(G')$ denote the set of all deterministic strategies for $G'$.

\paragraph{Statement of the adversarial game argument.}
A \emph{function on a Markov game} is any function on $\prod_{n\geq 0} \Omega_n$.
Given a morphism $\pi\colon G\to G'$, we may pull back a function $f$ on $G'$, to a function on $G$ whose value
at $(x_0,x_1,\dotsc)\in \prod_{n\geq 0}\Omega_n$ is $f(\pi_0(x_0),\pi_1(x_1),\dotsc)$. For ease of notation,
we shall denote the pull back of $f$ on $G$ by the same letter $f$. 

Though a general Markov game is not a probability space, we may naturally regard a Markov chain as a probability space.
So, for a Markov chain $M$, we write $\Pr_M$ and $\E_M$ for the probability and expectation on the probability space associated
to~$M$. In particular if $M=(\Omega,P)$ is a Markov chain, $P_n$ is the transition matrix for the $n$'th transition step
with respect to the probability~$\Pr_M$.

With these definitions at our disposal, we can formally state the idea at the start of this section.
\begin{theorem}[Adversarial game argument]\label{lem:adversarial}
Let $M\to G$ be a morphism from a Markov chain $M$ to a Markov game $G$. Suppose $f$ is any function on $G$. Then
\[
  \E_M [f]\geq \inf_{S\in \DStrat(G)} \E_S [f].
\]
\end{theorem}
The theorem has several uses beside bounding expectations from below:
\begin{itemize}
\item We can get an upper bound on $\E_M f$ by applying the theorem to $-f$.
\item  We may bound probabilities of suitable events in $M$ by applying the theorem to their characteristic functions.
\item By applying the theorem to the random variable $(f-\E_M f)^2$, and then using the inequality
$\E_S[ (f-c)^2 ]\geq \E_S[ (f-\E_S[f])^2 ]$ we also obtain
\[
  \Var_M[f]\geq \inf_{S\in \DStrat(G)} \Var_S [f].
\]
\item  Similarly, by applying the theorem to the function $-\exp(tf)$ we may bound the moment-generating
function of $f$. This gives a way to prove strong tail bounds on $f$. 
\end{itemize}

\paragraph{Penalties.} In applications of the adversarial game argument, we will sometimes have
constraints on adversary's action that are cumbersome to enforce using sets $R_n$. In this situation,
it will be convenient to impose a penalty instead. Formally, a \emph{penalty} is a function on a Markov game 
that is equal to some large constant when some condition is violated, and is equal to $0$ otherwise.
When $f$ is a function that the adversary tries to minimize and $\Penalty$ is a penalty, 
we instead consider a new function $f+\Penalty$. Similarly, if the adversary tries to maximize $f$, then the
new function is $f-\Penalty$. 

We shall always denote the penalty by $\Penalty$.

\paragraph{Compositions of morphisms.} The morphisms can be composed, i.e., whenever $\pi\colon G\to G'$ and $\pi'\colon G'\to G''$ are morphisms from $G$ to $G'$ and from $G'$ to $G''$ respectively, 
their composition $\pi'\circ\pi$ is a morphism from $G$ to $G''$. We omit the routine proof. 

\begin{observation}\label{obs:homimage} If $G\to G'$ is a morphism of Markov games and $f$ is a function on $G'$, then
$\inf_{S\in \DStrat(G)} \E_S [f]\geq \inf_{S'\in \DStrat(G')} \E_{S'} [f]$.
\end{observation}
Hence, for the purposes of the adversarial game argument, we may replace any game by its homomorphic image.
\begin{proof}[Proof of \Cref{obs:homimage}]
Let $S\in\DStrat(G)$ be arbitrary. Composition of morphisms $S\to G$ and $G\to G'$ yields a morphism $S\to G'$.
By the adversarial game argument $\E_S[f]\geq \inf_{S'\in \DStrat(G')} \E_{S'} [f]$. As $S$ is arbitrary, we are done.
\end{proof}

\paragraph{Proof of the adversarial game argument.}
We shall break the proof of \Cref{lem:adversarial} into two steps. We begin by showing a weaker version
of the lemma allowing randomness in adversary's strategy. Then we use derandomization
to obtain the full strength of the lemma.

\begin{lemma}[Adversarial game argument, basic version]
Let $M\to G$ be a morphism from a Markov chain $M$ to a Markov game $G$. Suppose $f$ is any function on $G$. Then
\[
  \E_M [f]\geq \inf_{S\in \Strat(G)} \E_S [f].
\]
\end{lemma}
\begin{proof}
Write $M=(\Omega,P)$ and $G=(\Adv',\Omega',P',R')$ and denote by $\pi\colon M\to G$ the morphism from
$M$ to~$G$. To define strategy $S=(\Omega',\widetilde{P})$ for $G$, we must specify the transition probabilities $\widetilde{P}_n$.
We put, for $x_n'\in \Omega_n'$ and $x_{n+1}'\in \Omega_{n+1}'$,
\begin{align*}
  \widetilde{P}_n(x_{n+1}',x_n')&\eqdef P_n'(x_{n+1}',x_n')\qquad &&\text{if }n\notin \Adv',\\
  \widetilde{P}_n(x_{n+1}',x_n')&\eqdef \Pr\nolimits_M[\pi_{n+1}(X_{n+1})=x_{n+1}'\mid \pi_n(X_n)=x_n'\bigr]\qquad&&\\
  &= \sum_{\substack{x_n\in \pi_n^{-1}(x_n')\\\mathclap{x_{n+1}\in \pi_{n+1}^{-1}(x_{n+1}')}}} P_n(x_{n+1},x_n)\Pr\nolimits_M\bigl[X_n=x_n\mid \pi_n(X_n)=x_n'\bigr]\qquad&&\text{if }n\in \Adv',
\end{align*}
where $(X_0,X_1,\dotsc,X_n,X_{n+1},\dotsc)$ is sampled from $\prod_{n\geq 0} \Omega_n$ using the Markov chain~$M$.

\subparagraph{$S$ is a Markov chain. }
We first must check that the functions $\widetilde{P}_n$ define a Markov chain on~$\Omega'$. To that end we need to verify that
$\sum_{x_{n+1}'\in \Omega_{n+1}'} \widetilde{P}_n(x_{n+1}',x_n')=1$ for each $x_n'\in\Omega_n'$. If $n\notin\Adv'$, this follows from $\widetilde{P}_n=\nobreak P_n'$ and
the respective property of $P_n'$. On the other hand, if $n\in\Adv'$, then this follows because the quantity $\Pr_M[\ \cdot \mid\nobreak \pi_n(X_n)=x_n']$ 
is a conditional probability on an actual probability space, and so the total probability is~$1$.


\subparagraph{$S$ is a valid strategy for $G$.}
We must verify that the identity maps on $(\Omega_n')_n$ induce a morphism
from $S=(\Omega',\widetilde{P})$ to the game $G=(\Adv',\Omega',P',R')$.  Since the condition \ref{morphism:inclusion} holds vacuously for morphisms from
Markov chains (because $\Adv=\emptyset$), we need only to 
check conditions \ref{morphism:mild} and \ref{morphism:strong}.\smallskip

We first check \ref{morphism:mild}. Suppose $n\in \Adv'$. Assume that the states $x_n'\in \Omega_n'$ and $x_{n+1}'\in \Omega_{n+1}'$ satisfy $\widetilde{P}_n(x_{n+1}',x_n')>0$.
From the definition of $\widetilde{P}_n$ in the case $n\in \Adv'$, it follows that there exist $x_n\in \pi_n^{-1}(x_n')$ and
$x_{n+1}\in \pi_n^{-1}(x_n')$ such that $P_n(x_{n+1},x_n)\Pr_M\bigl(X_n=x_n\mid \pi_n(x_n)=x_n'\bigr)>0$.
In particular, $P_n(x_{n+1},x_n)>0$. By the condition \ref{morphism:mild} applied to the morphism $\pi\colon M\to G$, it follows
that $\pi_{n+1}(x_{n+1})\in R_n'(\pi_n(x_n))$. As $\pi_n(x_n)=x_n'$ and $\pi_{n+1}(x_{n+1})=x_{n+1}'$, this implies that
$x_{n+1}'\in R_n'(x_n')$.

The condition \ref{morphism:strong} is nearly trivial to check. Because the morphism $S\to G$ is used by the identity maps on $(\Omega_n)_n$,
the condition \ref{morphism:strong} reduces to assertion that $\widetilde{P}_n=P_n'$ for $n\notin \Adv'$.

\subparagraph{Checking that $\E_M f=\E_S f$.} Use Markov chain $M$ to sample sequence $X=(X_0,X_1,\dotsc)$ from $\prod_{n\geq 0} \Omega_n$,
and use $S$ to sample $X'=(X_0',X_1',\dotsc)$ from $\prod_{n\geq 0}\Omega_n'$. Let $x'\in \prod_{n\geq 0} \Omega_n'$
be arbitrary. Write $x'_{\leq n}$ for $(x_0',\dotsc,x_n')$. Define the notations $X_{\leq n}$ and $X_{\leq n}'$ similarly.
We shall write $\pi(x_{\leq n})$ for $\bigl(\pi_0(x_0),\pi_1(x_1),\dotsc,\pi_n(x_n)\bigr)$.
We will show that, for every $n$ and every choice of $x'$, 
\begin{equation}\label{eq:strategyequality}
  \Pr\nolimits_M [\pi(X_{\leq n})=x_{\leq n}']=\Pr\nolimits_S [X_{\leq n}'=x_{\leq n}'].
\end{equation}
This will clearly imply $\E_M f=\E_S f$.

We use induction on $n$. The base case $n=0$ holds because both $\Omega_0$ and $\Omega_0'$ are single-element sets. Assume that
\eqref{eq:strategyequality} has been shown for all numbers that are at most~$n$.
Using the Markov property of $S$ we compute
\begin{align*}
  \Pr\nolimits_M [\pi(X_{\leq n+1})=x_{\leq n+1}']&=\Pr\nolimits_M[\pi_{n+1}(X_{n+1})=x_{n+1}' \mid \pi(X_{\leq n})=x_{\leq n}']\Pr\nolimits_M[\pi(X_{\leq n})=x_{\leq n}']\\
                   &=\Pr\nolimits_M[\pi_{n+1}(X_{n+1})=x_{n+1}' \mid \pi_n(X_n)=x_{n}']\Pr\nolimits_S[X'_{\leq n}=x_{\leq n}'].\\
\intertext{To complete the proof it suffices to show that $\Pr\nolimits_M[\pi_{n+1}(X_{n+1})=x_{n+1}' \mid \pi_n(X_n)=x_{n}']=\widetilde{P}_n(x_{n+1}',x_n')$, for it would then follow that}
\Pr\nolimits_S [\pi(X_{\leq n+1})=x_{\leq n+1}']&=\widetilde{P}_n(x_{n+1}',x_n')\Pr\nolimits_S[X'_{\leq n}=x_{\leq n}']&&\hskip-25ex\text{by Markov property of }S\\
&=\widetilde{P}_n(x_{n+1}',x_n')\Pr\nolimits_M[X'_{\leq n}=x_{\leq n}']&&\hskip-25ex\text{by the induction hypothesis}\\
                                  &=\Pr\nolimits_M[X_{\leq n+1}'=x_{\leq n+1}'].\\
\end{align*}

If $n\in \Adv'$, the requisite formula for $\Pr_M[\pi_{n+1}(X_{n+1})=x_{n+1}' \mid \pi_n(X_n)=x_{n}']$ follows from the definition of $\widetilde{P}_n(x_{n+1}',x_n')$.

If $n\notin \Adv'$, then we use the condition \ref{morphism:strong} for the morphism $\pi\colon M\to G$ to conclude
\begin{align*}
  \Pr\nolimits_M[&\pi_{n+1}(X_{n+1})=x_{n+1}' \mid \pi_n(X_n)=x_n']\\
         &=\sum_{x_n\in \pi_n^{-1}(x_n')} \Pr\nolimits_M[\pi_{n+1}(X_{n+1})=x_{n+1}' \mid X_n=x_n]\Pr\nolimits_M[X_n=x_n\mid \pi_n(X_n)=x_n'] \\
         &=\sum_{\substack{x_n\in \pi_n^{-1}(x_n')\\\mathclap{x_{n+1}\in \pi_{n+1}^{-1}(x_{n+1}')}}} P_n(x_{n+1},x_n)\Pr\nolimits_M[X_n=x_n\mid \pi_n(X_n)=x_n']\\
         &=\sum_{x_n\in \pi_n^{-1}(x_n')}P_n'\bigl(x_{n+1}',x_n'\bigr)\Pr\nolimits_M[X_n=x_n\mid \pi_n(X_n)=x_n']&&\!\!\!\!\!\!\text{by \ref{morphism:strong} for }\pi\colon M\to G\\
         &=\sum_{x_n\in \pi_n^{-1}(x_n')}\widetilde{P}_n\bigl(x_{n+1}',x_n'\bigr)\Pr\nolimits_M[X_n=x_n\mid \pi_n(X_n)=x_n']&&\!\!\!\!\!\!\text{since }n\notin \Adv'\\
         &=\widetilde{P}_n\bigl(x_{n+1}',x_n'\bigr). &&\qedhere
\end{align*}
\end{proof}

The second ingredient in the proof of \Cref{lem:adversarial} is a derandomization argument that permits
us to turn any randomized strategy into a deterministic strategy.

\begin{lemma}
Suppose $f$ is any function on a Markov game $G$, and $S$ is a strategy for $G$. Then there is a deterministic strategy
$S^*$ for $G$ satisfying
$
  \E_{S} [f]\geq \E_{S^*} [f].
$
\end{lemma}
\begin{proof}
Let $G=(\Adv,\Omega,P,R)$, and $S=(\Omega,\widetilde{P})$. Make a new probability space $\mathcal{M}$ consisting of
countably many independent random variables $\{D_{n,x_n} : n\in \Adv, x_n\in \Omega_n\}$.
The random variable $D_{n,x_n}$ is defined by $\Pr_{\mathcal{M}} [D_{n,x_n}=x_{n+1}]\eqdef \widetilde{P}_n(x_{n+1},x_n)$.
Let $\mathcal{S^*}=(\Omega,\widetilde{P}^*)$ be a deterministic strategy with the transition function
\[
\widetilde{P}^*_n(x_{n+1},x_n)=\begin{cases}
  1&\text{if }D_{n,x_n}=x_{n+1},\\
  0&\text{otherwise}.
\end{cases}
\]
Note that the strategy $\mathcal{S^*}$ is a random variable on the space $\mathcal{M}$.

The $\mathcal{S^*}$ is indeed a strategy for $G$, for the conditions \ref{morphism:mild} and \ref{morphism:strong} follow
from respective conditions for $\mathcal{S}$. Since
\[
  \E_{\mathcal{M}} \bigl[\E_{S^*} [f]\bigr]=\E_{S} [f],
\]
there is a choice of $S^*$ such that $\E_S [f]\geq \E_{S^*}[f]$.
\end{proof}

\section{Crude upper bound on \texorpdfstring{$\E[P_d-P_0]$}{E[P\textunderscore d-P\textunderscore 0]}}\label{sec:crudedriftupper}
\paragraph{Markov game strengthening.}
In this section we will prove a strengthening of \Cref{lem:crudedriftupper}. Its two
advantages are that it is easier
to prove, and that it is in the form that will be useful in the proof of the asymptotics 
for $\E[P_d-P_0]$ in \Cref{sec:asympt}.

The strengthening concerns the following Markov game.
\begin{game}\label{game:particles}
\gamestate Nonnegative integers $P_0<P_1<\dotsc<P_{\hl}$, which are interpreted as particle positions.
\gamestart $P_i=i$ for all $i=0,1,\dotsc,\hl$.
\gameduration{$L$}
\gameturn Adversary chooses a partition $\mathcal{A}$ of $\{0,1,\dotsc,\hl\}$ into $k$ parts, some of which are possibly empty.
Then Fortune picks a set $A\in\mathcal{A}$ uniformly among the $k$ sets in $\mathcal{A}$. The particle positions are
then updated as in the \Cref{step:Aadvance} of \Cref{prop:update}.
\end{game}
Since this is our first example of a Markov game, before stating our result we explain how to translate this description
into the formalism of the preceding section. 

We use the personified `Fortune' to denote the times when the transitions
are random, and `Adversary' to denote the remaining times. In this game, actions of Adversary and Fortune alternate,
with Adversary acting first. 
Hence, $\Adv=\{0,2,4,6,\dotsc\}$ is the set of times when Adversary acts. 
Since each game turn is made of two actions, and the game lasts $L$ turns, this means that we are 
interested only in the times $0,1,\dotsc,2L$, and can ignore whatever happens afterward.
For the same reason, the notation $P_i(L)$ below refers to the value of $P_i$ at the end of the $L$'th turn, i.e.,
at the time~$2L$.

At the beginning of Adversary's turns, the state space consists of an ordered tuple $(P_0,P_1,\dotsc,P_{\hl})$ satisfying $P_0<P_1<\dotsb<P_d$.
After Adversary makes their turn, the state space becomes enriched with the partition $\mathcal{A}$ that they choose.
Hence, denoting by $\mathcal{P}(\hl)$ the set of of all partitions of $\{0,1,\dotsc,\hl\}$ into exactly $k$ (possibly empty) parts, we
have
\begin{align*}
  \Omega_0&=\{(0,1,\dotsc,\hl)\}\\
  \Omega_n&=\{(P_0,P_1,\dotsc,P_{\hl})\in \Z_{\geq 0}^{\hl} : P_0<P_1<\dotsc<P_d \}&&\text{if }n\in \Adv\setminus\{0\},\\
  \Omega_n&=\{(P_0,P_1,\dotsc,P_{\hl},\mathcal{A})\in \Z_{\geq 0}^{\hl}\times \mathcal{P}(\hl) : P_0<P_1<\dotsc<P_d \}&&\text{if }n\notin \Adv.
\end{align*}
The functions $R_n$ are simple: $R_n(P_0,P_1,\dotsc,P_{\hl})=\{(P_0,P_1,\dotsc,P_{\hl})\}\times \mathcal{P}(\hl)$.
The definition of $P_n$ is more cumbersome: suppose $n\notin\Adv$,
and let $x_n=(P_0,P_1,\dotsc,P_{\hl},\mathcal{A})\in \Omega_n$.
For a set $A\in\nobreak\mathcal{A}$, denote by $\operatorname{Evolve}(x_n,A)$ the
particle positions after applying the procedure from the \Cref{step:Aadvance} of \Cref{prop:update}
to particles in positions $P_0,P_1,\dotsc,P_{\hl}$ with the set $A$ in place of $A[L]$.
For a state $x_{n+1}=(P_0',P_1',\dotsc,P_{\hl}')\in \nobreak\Omega_{n+1}$, let $r$ be the number of sets
$A\in\mathcal{A}$ for which $\operatorname{Evolve}(x_n,A)=x_{n+1}$. Then $P_n(x_{n+1},x_n)=r/k$.\medskip

%
%
%
%
%


We are now ready to state the main result of this section.
\begin{lemma}\label{lem:adversariallndsupper}
For every adversary's strategy in \Cref{game:particles}, $\E[P_{i+1}(L)-P_i(L)]\leq \sqrt{2L/k}+1$
for all~$i$. In particular, $\E[P_{\hl}(L)-P_0(L)]\leq \hl\sqrt{2L/k}+\hl$ for every strategy.
\end{lemma}
The expectation here is taken with respect to the strategy. In particular,
the lemma immediately implies \Cref{lem:crudedriftupper}. Indeed, let $G=(I,\Omega,P,R)$ be the formalization of \Cref{game:particles} 
constructed above, and consider the strategy $S=(\Omega',P')$ where 
$\Omega_n'$ consists of tuples $(P_0,P_1,\dotsc,P_{\hl},w,w')$ where $P_0<P_1<\dotsb<P_{\hl}$ are nonnegative integers and 
$w,w'$ are words over $[k]$ of lengths $P_{\hl}+1$ and $\lfloor n/2\rfloor$ respectively. The transition
from $\Omega_{2\ell}'$ to $\Omega_{2\ell+1}'$ is the identity, i.e., the chain remains in the same state.
The transition from $\Omega_{2\ell+1}'$ to $\Omega_{2(\ell+1)}'$ is performed by appending a random symbol to $w'$ and then executing the algorithm in 
\Cref{prop:update} with $\ell$ in place of~$L$. This is indeed a strategy for \Cref{game:particles}, with the morphism $\pi$
given by $\pi_{2\ell}(P_0,P_1,\dotsc,P_{\hl},w,w')=(P_0,P_1,\dotsc,P_{\hl})$ for even times and by 
$\pi_{2\ell+1}(P_0,P_1,\dotsc,P_{\hl},w,w')=(P_0,P_1,\dotsc,P_{\hl},\mathcal{A}[\ell])$ for odd times.

We have stated the lemma in terms
of $\hl+1$ particles $P_0,\dotsc,P_{\hl}$ rather than the pair $P_i$, $P_{i+1}$ only because we wanted to
point its implication on the quantity $\E[P_{\hl}(L)-P_0(L)]$. For the purpose of giving a proof, we may however strip away the irrelevant
particles, and focus solely on the gap between a fixed pair of consecutive particles.
We obtain a much simpler game:

\begin{game}\label{game:delta}
\gamestate Positive integer $\Delta$ (which we interpret as a gap between two consecutive particles).
\gamestart $\Delta=1$.
\gameduration{$L$}
\gameturn 
\begin{itemize}[leftmargin=*]
\item First, Adversary chooses a vector $v$ that either is the zero vector $\vec{0}\in\Z^k$ or is a permutation
of the vector $(+1,-1,0,0,0,\dotsc,0)\in\Z^k$. 
\item \label{step:decrement} Second, Adversary decrements one of the coordinates of $v$, or keeps it intact.
\item Third, Fortune then picks $i\in [k]$ uniformly at random, and adds $v_i$ to $\Delta$. 
\item Finally, if $\Delta$ is $0$ or less, then we reset $\Delta$ to $1$.
\end{itemize}
\gameobjective Adversary aims to maximize $\Delta$.
\end{game}
In the formal language of \Cref{sec:adversarial}, we may say that \Cref{game:particles} admits a morphism
into \Cref{game:delta}. Again, we spell out the details.

We model \Cref{game:delta} similarly to how we modeled \Cref{game:particles}:
it is the Markov game $(\Adv',\Omega',P',R')$, where $\Adv'=\{0,2,4,6,\dotsc\}$ (same
as in \Cref{game:particles}), $\Omega_n'=\N$ for positive even $n$, whereas for odd values of $n$ the set $\Omega_n'$ consists of pairs 
of the form $(\Delta,v)$ where $\Delta\in\N$ and $v\in V$ with 
\begin{align*}
V&\eqdef U+W,\\
U&\eqdef\{\vec{0}\}\cup\{v\in \Z^k: v \text{ is a permutation of }(+1,-1,0,0,0,\dotsc,0)\},\\
W&\eqdef\{\vec{0}\}\cup\{v\in \Z^k: v \text{ is a permutation of }(-1,0,0,0,0,\dotsc,0)\}.
\end{align*}
When $n$ is odd, $(\Delta,v)\in \Omega_n'$, $\Delta'\in \Omega_{n+1}'$,
we set $P_n'\bigl(\Delta',(\Delta,v)\bigr)=r/k$ where $r$ is the number of coordinates
$i\in [k]$ such that $\Delta'=\max(\Delta+v_i,1)$. Finally,
the $R_n'$ is defined by $R_n'(\Delta)=\{\Delta\}\times V$ (for even $n$).

The morphism $\pi$ from \Cref{game:particles} to \Cref{game:delta}
maps the $(\hl+1)$-tuple $(P_0,\dotsc,P_{\hl})$ into the number $\Delta=P_{i+1}-P_i$ and the partition $\mathcal{A}$
into vector $v$ in such a way that choice of $i$'th part in $\mathcal{A}$ corresponds to
choice of $v_i$. The decrementation in \cref{step:decrement} corresponds to the particle
$P_{i-1}$ bumping into $P_i$ and pushing $P_i$.

\paragraph{Proof of \texorpdfstring{\Cref{lem:adversariallndsupper}}{Lemma \ref{lem:adversariallndsupper}}.}
Invoking \Cref{obs:homimage}, we see that \Cref{lem:adversariallndsupper} would follow
if we show that $\E[\Delta]\leq 1+\sqrt{2L/k}$ for every strategy in \Cref{game:delta}. That is
precisely what we will do.

Let $\Delta'\eqdef \Delta-1/2$. It suffices to show that 
\begin{equation}\label{eq:deltaprime}
  \E[\Delta'(L)^2]\leq \tfrac{1}{4}+2L/k
\end{equation}
holds for every Adversary's strategy. Indeed, 
Cauchy--Schwarz inequality would then imply that
$\E[\Delta(L)]\leq \sqrt{\tfrac{1}{4}+2L/k}+\tfrac{1}{2}\leq \sqrt{2L/k}+1$.

We prove \eqref{eq:deltaprime} by induction on $L$. The base case $L=0$ is immediate, so assume $L>0$. 

Consider the vector $v$ chosen by Adversary at the last (i.e., the $L$'th) turn
of the game. Write it as $v=u+w$ where $u$ is either $\vec{0}$ or a permutation of $(+1,-1,0,0,\dotsc,0)$,
and $w$ is either $\vec{0}$ or a permutation of $(-1,0,\dotsc,0)$.
Fix any strategy $S$ for \Cref{game:delta}, and denote by $i\in [k]$ the Fortune's random choice at turn $L$.
We compute
\begin{align*}
  \E_S\bigl[(\Delta'(L-1)+u_i)^2\bigr]&\leq \sum_{x\in \Z+1/2} \Pr\nolimits_S[\Delta'(L-1)=x]\bigl(\tfrac{1}{k} (x-1)^2+\tfrac{1}{k}(x+1)^2+\tfrac{k-2}{k}x^2\bigr)\\
                                  &=\sum_{x\in \Z+1/2}\Pr\nolimits_S[\Delta'(L-1)=x](x^2+\tfrac{2}{k})\\
                                  &=\E_S\bigl[\Delta'(L-1)^2\bigr]+\tfrac{2}{k}.
\end{align*}
Since $\Delta'(L)=\max(\Delta'(L-1)+u_i+w_i,1/2)\leq \abs{\Delta'(L-1)+u_i}$, it thus follows that
\[
  \E_S\bigl[\Delta'(L)^2\bigr]\leq \E_S\bigl[\Delta'(L-1)^2\bigr]+\tfrac{2}{k},
\]
concluding the induction step.

\section{Crude lower bound on \texorpdfstring{$\E[P_1-P_0]$}{E[P\textunderscore 1-P\textunderscore 0]}}\label{sec:crudedriftlower}
Similarly to the proof of the upper bound in the preceding section we shall compare
the evolution of $P_1-P_0$ to an (adversarial) random walk. The appropriate random
walk is lazy, i.e., there is a non-negligible probability that in a given step nothing
happens. To argue that $\E[P_1-P_0]$ is large, we must show that the random walk is not
too lazy. This makes the argument more complex compared to that in the previous section.

The game we use is similar to the two-particle version of \Cref{game:particles}, except that
we must keep track of whether the symbols $w[P_0]$ and $w[P_1]$ are same or different.

\begin{game}\label{game:two}
\gamestate Pair $(\Delta,F)$ where $\Delta\in \N$ and $F\in \{\same,\diff\}$.
\gamestart $\Delta=1$, whereas 
$F=\same$ with probability $1/k$ and $F=\diff$ with probability $1-1/k$.
\gameduration{$L$}
\gameturn 
What happens depends on the value of $F$ at the turn's start.
\begin{itemize}[leftmargin=*]
\item If $F=\same$, then with probability $1/k$ Fortune decides to toss a coin.
\item If $F=\diff$, then 
\begin{itemize}
\item with probability $1/k$, Fortune decrements $\Delta$ (setting it to $1$ afterward should it become $0$), and Adversary chooses the new value of $F$,
\item with probability $1/k$, Fortune increments $\Delta$, and decides to toss a coin,
\item with probability $1-2/k$, nothing happens.
\end{itemize}
\item The coin, should Fortune decide to toss it, lands on heads with probability $1-1/k$.
Should it land on heads, the new value of $F$ is set to $\diff$, and
should it land on tails, the new value of $F$ is set to $\same$.
\end{itemize}
\gameobjective Adversary aims to minimize $\Delta$.
\end{game}
There is a morphism from the Markov chain of \Cref{prop:update} into this game, which can informally be described as setting
$\Delta=P_1-P_0$ and setting $F=\same$ iff $w[P_0]=w[P_1]$.

The formal treatment is similar to that of \Cref{game:particles,game:delta}
in the preceding section: each turn of \Cref{game:two} is represented by two time steps, the first step being random, and the second step
being adversarial. The exception is the very first step in which the initial value of $F$ is chosen. So,
$\Adv=\{2,4,6,\dotsc\}$. 
The state spaces are 
\begin{align*}
  \Omega_n&=\{(\Delta,F) : \Delta\in \N,\ F\in \{\same,\diff\}\}&&\text{if }n\notin \Adv\text{ and }n\geq 1,\\
  \Omega_n&=\{(\Delta,F,C) : \Delta\in \N,\ F\in \{\same,\diff\},\ C\in\{\choice,\nochoice\} \}&&\text{if }n\in \Adv,
\end{align*}
where the variable $C$ records whether the Adversary has the choice of $F$ at the next step.
The functions $R_n$ are defined by 
\begin{align*}
R_n(\Delta,F,\nochoice)&=\{(\Delta,F)\},\\
R_n(\Delta,F,\choice)&=\{(\Delta,\same),(\Delta,\diff)\}.
\end{align*}
The transition functions $P_n$ are defined by
\begin{align*}
  P_n\bigl( (\Delta,\diff,\nochoice) , (\Delta,\same) \bigr)&=\tfrac{1}{k}(1-\tfrac{1}{k}),\\
  P_n\bigl( (\Delta,\same,\nochoice) , (\Delta,\same) \bigr)&=1-\tfrac{1}{k}+\tfrac{1}{k^2},\\
  P_n\bigl( (\max(\Delta-1,1),\same,\choice) , (\Delta,\diff) \bigr)&=\tfrac{1}{k},\\
  P_n\bigl( (\Delta+1,\diff,\nochoice) , (\Delta,\diff) \bigr)&=\tfrac{1}{k}(1-\tfrac{1}{k}),\\
  P_n\bigl( (\Delta+1,\same,\nochoice) , (\Delta,\diff) \bigr)&=\tfrac{1}{k^2},\\
  P_n\bigl( \Delta,\diff,\nochoice) , (\Delta,\diff) \bigr)&=1-\tfrac{2}{k}.
\end{align*}

Call a turn in which Fortune either decrements or increments $\Delta$ a \emph{good} turn.
Let $H$ be the total number of coin tosses that landed on heads. Let
$G$ be the number of good turns. Once a coin lands on heads, 
at least one good turn must occur before the next coin toss.
Hence, $G\geq H-1$.

Because each turn Fortune decides to toss a coin with the same probability $1/k$,
it follows that $H\sim\nobreak \Binom\bigl(L,(1/k)(1-1/k)\bigr)$. Since $(1/k)(1-1/k)\geq 1/2k$ (in view of $k\geq 2$), 
from an asymmetric version of the Chernoff bound (see \cite[Theorem~A.1.13]{alon_spencer}) we infer that
\begin{equation}\label{eq:asymchernoffL}
  \Pr[G\leq L/4k-1]\leq \Pr[H\leq L/4k] \leq \exp\Bigl(-\frac{(L/4k)^2}{2(1/k)(1-1/k)L}\Bigr)\leq \exp (-L/32k).
\end{equation}
Because of this estimate, we may simplify the game further:
\begin{game}\label{game:twosimple}
\gamestate Nonnegative integer $\Delta$. 
\gamestart $\Delta=1$. 
\gameduration{$L$}
\gameturn Adversary decides if they want this turn to be good. If they decide on the turn
being good, Fortune adds either $-1$ or $+1$ to $\Delta$ at random (setting
$\Delta$ to $1$ afterward should it become $0$). If they decide on the turn being bad,
nothing happens.
\gameobjective Adversary aims to minimize $\Delta$.
\gamepenalty Adversary pays penalty of $T+1$ to the objective function if the total number of good turns is less than $T\eqdef L/4k-1$.
\end{game}
Recall that penalties provide a shorthand for defining objective functions. So, formally the objective 
is equal to $\Delta$ if the number of good turns is at least $L/4k-1$, and is equal to $\Delta+T+1$ otherwise.
However, as we shall shortly see, separating the penalty is convenient in the analysis.

If Adversary uses an optimal strategy, then the penalty condition is never triggered. Indeed,
if $S$ is any strategy that sometimes triggers the penalty, we may modify it so it does not.
Namely, when exactly $t$ turns are left and only $T-t$ good turns have been played,
not allowing the next turn to be good triggers the penalty. In any such situation, Adversary
should make all the remaining turns good. Doing so strictly improves the strategy.
So, every optimal strategy always uses at least $T$ good turns.

Because it does not matter which turns are good, we may assume that the first $T$ turns are good,
whereas the subsequent turns are subject to Adversary's strategy. We claim that Adversary's best strategy
is to disallow good turns on all of these subsequent turns.

Indeed, fix any strategy and denote by $\Delta(\ell)$ the value of $\Delta$ at the end of the $\ell$'th turn.
Then $\Delta(\ell+1)=\max(\Delta(\ell)+v_\ell,1)$ where $v_{\ell}$ is either $0$ (if the turn is bad)
or a uniform element of $\{\pm 1\}$. In either case, $\E[\Delta(\ell+1)]\geq \E[\Delta(\ell)]$
with equality if the turn is bad. Hence, $\E[\Delta(L)]\geq E[\Delta(T)]$ with equality
if all but the first $T$ turns are bad. So, in particular, Adversary's optimal strategy is to 
disallow good turns.\smallskip

The reflection principle tells us that the value of $\Delta$ under the optimal play
is the same as \linebreak$1/2+\abs{\Delta'}$ where $\Delta'$ is the position of a simple random walk on
$\Z+1/2$ starting from~$1/2$ after $T$ steps. Therefore, 
\[
\E[\Delta]\geq \E\bigl[1/2+\abs{\Delta'}\bigr]=
\begin{dcases}
\frac{1}{2}+2^{-T}\cdot \Bigl(T+\frac{1}{2}\Bigr)\binom{T}{T/2},&\text{if }T\text{ is even},\\
\frac{1}{2}+2^{-T}\cdot 2T\binom{T-1}{(T-1)/2},&\text{if }T\text{ is odd},
\end{dcases}
\]
where the formula for $\E\bigl[\abs{\Delta'}\bigr]$ is proved in \Cref{sec:randomwalk}.
Both in the case $T$ is even and in the case $T$ is odd, using Stirling's formula with explicit error term (see \cite[Eq.~(9.91)]{graham_knuth_patashnik}), we may show that
$\binom{2n}{n}\geq e^{-1/6n}\frac{4^n}{\sqrt{\pi n}}$ for all $n$, and use this to deduce that
\[
  \E[\Delta]\geq \tfrac{1}{2}+\sqrt{\tfrac{4}{7}(T+1)} ,\qquad\text{if }L\geq 200k,
\]
in the optimal strategy for \Cref{game:twosimple}.

Let $\Penalty$ be the penalty function from \Cref{game:twosimple}. 
Recall that the pull back of $\Penalty$ from \Cref{game:twosimple} to \Cref{game:two} is a function
on \Cref{game:two}, which is also denoted by $\Penalty$.
Using \eqref{eq:asymchernoffL} and applying \Cref{obs:homimage} to the morphism from \Cref{game:two} to \Cref{game:twosimple}, we obtain
\begin{align*}
  \inf_{S\in\DStrat(\text{\Cref{game:two}})} \E_S[\Delta]+(T+1)\exp (-L/32k)&\geq \inf_{S\in\DStrat(\text{\Cref{game:two}})} \E_S[\Delta+\Penalty]\\
                                         &\geq \inf_{S'\in\DStrat(\text{\Cref{game:twosimple}})} \E_{S'}[\Delta+\Penalty]\\
                                         &\geq \frac{1}{2}+\sqrt{\frac{L}{7 k}}.
\end{align*}
Since $(T+1)\exp (-L/32k)\leq \tfrac{1}{2}$ when $L\geq 200k$, \Cref{lem:crudedriftlower} follows.

\section{Most expectant partitions are trivial}\label{sec:trivial}
As sketched in the introduction, the key to the proof of \Cref{thm:driftasympt} is to show that almost all expectant partitions are trivial,
which are defined as partitions all of whose non-empty parts are singletons.
Let $B\eqdef \{\ell< L : \cA[\ell]\text{ is non-trivial}\}$. This section is devoted to the proof of the following estimate.
\begin{lemma}\label{lem:trivialpartitionbound}
The number of non-trivial expectant partitions satisfies 
\[
  \Pr[\abs{B}\geq 6\hl^2L/k ]\leq 4\hl^2 L \exp(-L^{1/2}k^{-3/2})
\]
for all $\hl,L\geq 1$ and all $k\geq 2$.
\end{lemma}

\paragraph{Bumps and jumps.}
We shall define $Q_0,\dotsc,Q_{\hl}$ in such a way that $Q_0,\dotsc,Q_{\hl}$ is a permutation of $P_0,\dotsc,P_{\hl}$.
In particular, an expectant partition will be non-trivial if and only if $w[Q_i]=w[Q_j]$ for some pair $i\neq j$.

Following the interpretation of $P_0,\dotsc,P_{\hl}$ as particles, we shall think of $Q_0,\dotsc,Q_{\hl}$ also as particles. 
Informally, we think of $P$-particle that tries to move into an already-occupied positions as bumping the particle that
is located there. On the other hand, $Q$-particles will avoid moving into an already-occupied position, jumping over
the particles in front of it instead.

Formally, we begin with $Q_i(0)=i$ for $i=0,1,\dotsc,\hl$. At time step $\ell$, we examine $w'[\ell]$, and use its value to define the set $I=\{i : w[Q_i(\ell)]=w'[\ell]\}$.
We think of particles $\{Q_i : i\in I\}$ as being excited. We then examine excited particles in decreasing order
of their positions. We move each excited particle to next vacant position on the right, jumping over non-excited particles if necessary.

The advantage of $Q$-particles over $P$-particles is that we may focus on a single pair of particles at a time, without worrying that other particles might bump 
and displace them. 

The figure below illustrates the difference between the dynamics of $P$-particles and $Q$-particles, during the same time step.

\begin{center}
\def\ldist{0.45} 
\def\sdist{6}    
\def\gdist{0.18} 
\newcommand*{\partpic}[8]{%
\foreach \i in {0,1,...,9}
   \draw ($\i*(\ldist,0)+#2*(\sdist,0)$) circle(1.2pt);
\foreach \i / \j in {#3/0,#4/1,#5/2,#6/3,#7/4,#8/5}
{
  \StrLeft{\i}{1}[\xind]
  \pgfmathsetmacro{\xpos}{\xind*\ldist+#2*\sdist}
  \begin{scope}[xshift=\xpos cm]
  \IfStrEq{\i}{\xind}
  {\fill (0,0) circle (2.25pt);}
  { 
    \fill [scale=0.0035] (-21.0730,  12.6550) -- (-28.4610,  27.5920) -- (-10.1590,  18.2990) -- (-0.4280,  37.1280) -- ( 8.1530,  16.7350) -- ( 35.3860,  27.7750) -- ( 22.6170,  8.9130) -- ( 46.1020,  6.8930) -- ( 22.6590, -1.9590) -- ( 41.8300, -21.4980) -- ( 16.4580, -13.5130) -- ( 13.9380, -32.0150) -- ( 1.4620, -18.2710) -- (-8.4830, -36.5670) -- (-11.7510, -16.6220) -- (-33.6870, -27.3790) -- (-25.1500, -6.8760) -- (-44.1890, -6.4460) -- (-27.1700,  2.5050) -- (-46.1020,  16.6100) -- cycle;
  }
  \end{scope}
  \node at (\xpos,0.35) {$\scriptstyle #1_{\j}$};
}}
\newcommand*{\statepic}[8]{%
  \draw[very thin, rounded corners=2pt] ($(-0.4,-0.3) + #2*(\sdist,0)$) rectangle ($9*(\ldist,0)+(0.3,2.1)+#2*(\sdist,0)$);
  \node at ($(0.05,-0.7) + 4.5*(\ldist,0) + #2*(\sdist,0)$) {#1}; 
  \begin{scope}[yshift=1.5cm]
  \partpic{P}{#2}{#3}{#4}{#5}{#6}{#7}{#8}
  \end{scope}
  \partpic{Q}{#2}
}
\begin{tikzpicture}
\statepic{Previous state}{0}{0}{2}{3}{6}{7}{8}{3}{7}{0}{6}{2}{8}
\statepic{Four particles are excited}{1}{0x}{2x}{3}{6x}{7x}{8}{3}{7x}{0x}{6x}{2x}{8}
\statepic{Next state}{2}{1}{3}{4}{7}{8}{9}{3}{9}{1}{7}{4}{8}
\foreach \i/\j in {0/1,1/2}
   \draw[->] ($9*(\ldist,0)+(0.3,0.9)+(\gdist,0)+\i*(\sdist,0)$) -- ($(-0.4,0.9)+(-\gdist,0)+\j*(\sdist,0)$); 
\node at ($(0.05,-1.6) + 4.5*(\ldist,0) + 1*(\sdist,0)$) {\textbf{Figure 2: }Evolution of $P$- and $Q$-particles (example).};
\end{tikzpicture}
\end{center}

\paragraph{Two-particle evolution.}
Fix a pair $i\neq j$; and let $B_{ij}\eqdef \{\ell : w[Q_i(\ell)]=w[Q_j(\ell)]\}$. We will
show that
\begin{equation}\label{eq:trivialpairtail}
  \Pr[\abs{B_{ij}}\geq 6L/k]\leq 2L \exp(-L^{1/2}k^{-3/2}).
\end{equation}
Since $B=\bigcup_{i\neq j} B_{ij}$ and there are $\binom{\hl+1}{2}$ pairs $i\neq j$, \Cref{lem:trivialpartitionbound} will then follow from the union bound.

Throughout the rest of the section, we use notations $\Qmax\eqdef\max(Q_i,Q_j)$ and $\Qmin\eqdef \min(Q_i,Q_j)$.
Since the right side of \eqref{eq:trivialpairtail} exceeds $1$ if $L\leq k^3$, we
may also assume that $L>k^3$ in what follows.

Since $Q$-particles move if and only if they are excited, it is tempting
to immediately discard all the particles except for $Q_i$ and $Q_j$. This requires
a little care because $Q_i$ and $Q_j$ might still jump over the other particles.
However, as the next lemma shows, the values of $w[\Qmax]$ still behave as if
the other particles do not exist.

\begin{lemma}\label{lem:Qijindep}
Let $\ell$ be arbitrary. Then
\[
  \Pr\bigl[w[\Qmax(\ell)]=s\mid \bigl(\Qmax(\ell)>\Qmax(\ell-1)\bigr) \wedge \Hist(\ell) \bigr]=1/k\qquad\text{for all }s\in [k],
\]
where $\Hist(\ell)$ consists of values of $5$-tuples $(Q_i(t),Q_j(t),w[Q_i(t)],w[Q_j(t)],w'[t])$
for all $t<\ell$.
\end{lemma}
In particular, if we consider the times $\ell$ when $\Qmax(\ell)>\Qmax(\ell-1)$, and write down the
sequence of values of $w[\Qmax(\ell)]$ at those times, then we obtain a uniform random word. 

\begin{proof}[Proof of~\Cref{lem:Qijindep}]
It suffices to prove, for each $T$, that
\[
  \Pr\bigl[w[\Qmax(\ell)]=s\mid \bigl(T=\Qmax(\ell)>\Qmax(\ell-1)\bigr) \wedge \Hist(\ell) \bigr]=1/k\qquad\text{for all }s\in [k].
\]

Consider \emph{$T$-bounded} dynamics, which is identical to the $Q$-particle dynamics except that
once the value of $Q_r$, for $r\neq i,j$, becomes $T$ or larger, we stop tracking~$Q_r$. We 
imagine uncovering the symbols of $w$ only when a particle lands on it. Under this coupling,
if $\ell$ is the first time when $\Qmax$ becomes $T$ or larger, the value of $w[\Qmax(\ell)]$ 
is not revealed before time~$\ell$, and so is independent of the history before the time $\ell$.
\end{proof}

We can delay exposing $w[\Qmax(\ell)]$ for even longer than in the preceding proof:
Imagine that for each symbol $w[p]$ of $w$ we generate a random $\{0,1\}$-vector of length $k$
having a single $1$ and $k-1$ many $0$'s. We then set $w[p]$ to the position of the single $1$
in that vector. In the $T$-bounded dynamics, when $\Qmax$ becomes equal to some value $T$, we do not expose $w[T]$ completely, 
but instead only check if $w[T]=w[\Qmin]$ by exposing $w[\Qmin]$'th position of the vector
associated to $w[T]$. That is enough to decide if $i$ and $j$ belong to the same part
of the expectant partition. Then each time $\Qmin$ increases, we similarly check for $w[\Qmax]=w[\Qmin]$.
This gives us a way to bound the conditional probability
\[
\Pr\bigl[ w[\Qmax(\ell)]=w[\Qmin(\ell)] \mid (\Qmin(\ell)>\Qmin(\ell-1)) \wedge \Hist(\ell)\bigr].
\]
Indeed, if $T=\Qmax(\ell)=\Qmax(\ell-1)$, then the probability is $0$
if $w[\Qmin(\ell)]$'th position of the vector associated to $w[T]$ has been exposed,
and otherwise the probability is $1/r$ where $r$ is the number of yet-unexposed entries in the vector associated to $w[\Qmax]$.
If $T=\Qmax(\ell)>\Qmax(\ell-1)$, then the probability is $1/k$. So, in all three cases, we obtain that
\[
 \Pr\bigl[ w[\Qmax(\ell)]=w[\Qmin(\ell)] \mid (\Qmin(\ell)>\Qmin(\ell-1)) \wedge \Hist(\ell)\bigr]\leq 1/r.
\]

This suggests the following Markov game, in which we give Adversary the power to choose new values of $Q_i$, $Q_j$ and $w[\Qmin]$, whenever
$Q_i$ and $Q_j$ change. Furthermore, when Adversary chooses the value of $w[\Qmin]$, they have the extra power of overwriting the old value of $w[\Qmin]$.
Adversary can also choose the initial values of $Q_i$ and $Q_j$. To aid the analysis, the initial set
of yet-unexposed positions in $w[\Qmax]$ will be arbitrary. By symmetry, only its size will matter.

More specifically, the game depends on two parameters $s_0\in \{0,1,\dotsc,k\}$ and 
$b_0\in\{\stin,\stout\}$ determining the initial state. In this game, $S$ corresponds to the set of yet-unexposed positions of the $\{0,1\}$-vector
associated to $w[\Qmax]$. The parameters $s_0$ and $b_0$ indicate the initial size of $S$, and whether $w[\Qmin]$ is an element of~$S$.
Though we are interested only in the case $(s_0,b_0)=(k,\stin)$, the extra generality will be useful in the analysis.

\begin{game}\label{game:trivial}
\gamestate Set $S\subseteq [k]$, two non-negative distinct integers $Q_i,Q_j$ and a symbol $w[\Qmin]\in [k]$.
\gameduration{$T$}
\gamestart Set $S=[s_0]$ (if $s_0=0$, then $S=\emptyset$). Adversary chooses the initial values of 
$Q_i$ and $Q_j$. Adversary also chooses the initial value of $w[\Qmin]$, which has to be an element of $S$ if $b_0=\stin$,
or an element of $[k]\setminus S$ if $b_0=\stout$. 
\gameturn
\begin{enumerate}[label=(G\arabic*), ref=(G\arabic*),leftmargin=*]
\item First, if $w[\Qmin]\in S$, then 
  \begin{itemize}
  \item with probability $1/\abs{S}$, Fortune sets $S$ to $\emptyset$,
  \item with probability $1-1/\abs{S}$, Fortune removes $w[\Qmin]$ from the set $S$.
  \end{itemize}
  If, on the other hand, $w[\Qmin]\notin S$, nothing happens.
\item Second,\label{gamestep:main}
\begin{enumerate}[label=(\alph*), ref=(G\arabic{enumi}\alph*)]
  \item \label{gamesubstep:one} if $S\neq \emptyset$, then 
    \begin{itemize}
    \item with probability $1/k$, Fortune forces Adversary to increase $Q_i$ (to a~value of Adversary's~choice),
    \item with probability $1/k$, Fortune forces Adversary to increase $Q_j$ (to a~value of Adversary's~choice),
    \item with probability $1-2/k$, nothing happens,
    \end{itemize} 
  \item \label{gamesubstep:two} if $S=\emptyset$, then
    \begin{itemize}
      \item with probability $1/k$, Fortune forces Adversary to increase both $Q_i$ and $Q_j$ to values of Adversary's choice,
      \item with probability $1-1/k$, nothing happens,
    \end{itemize}
  \item if $\Qmax$ is increased as a result of \ref{gamesubstep:one} or \ref{gamesubstep:two}, then $S$ is reset to $[k]$,
  \item if $\Qmin$ or $\Qmax$ is increased as a result of \ref{gamesubstep:one} or \ref{gamesubstep:two}, then Adversary sets $w[\Qmin]$ at will.
\end{enumerate}
\end{enumerate}
\gameobjective Let $\overline{B}(s_0,b_0,T)$ be the number of those turns for which
we have $S=\emptyset$ at the start of \cref{gamestep:main}. The Adversary's goal is to maximize $f=f\bigl(\overline{B}(s_0,b_0,T)\bigr)$, which is some non-decreasing function of $\overline{B}(s_0,b_0,T)$.
\end{game}

From the preceding discussion and the adversarial game argument, we know that
\begin{equation}\label{eq:trivialstrategyarg}
  \E[f]\leq \sup_{\mathcal{S}\in \DStrat(\text{\Cref{game:trivial}}
    )} \E_{\mathcal{S}}[f].
\end{equation}
In our application of \Cref{game:trivial}, $f$ will be the characteristic function of the event $\overline{B}(s_0,b_0,T)\geq 3pL$ for a suitable choice of $p\in [0,1]$.
For simplicity and generality of the analysis, we shall assume that $f$ is an arbitrary non-decreasing function of $\overline{B}(s_0,b_0,T)$.
Note that any such function is a function on \Cref{game:trivial} in the sense of \Cref{sec:adversarial}.

\paragraph{Analysis of \texorpdfstring{\Cref{game:trivial}}{Game \ref{game:trivial}}.}
Because the function $f$ that Adversary tries to maximize is non-decreasing, Adversary's optimal strategy
must be memoryless, i.e., Adversary's actions should depend only on the current set $S$, integers $Q_i,Q_j$ and the symbol $w[\Qmin]$.
Furthermore, from the symmetry, it is clear that only the size of $S$ matters, not the actual constituent elements.
It is not hard to see, by induction on $T$, that decreasing the size of $S$ is advantageous to Adversary,
i.e., decreasing $\abs{S}$ can only increase the value of $\overline{B}(s_0,b_0,T)$. From this it follows that
Adversary should set $w[\Qmin]$ to an element of $S$ whenever they can. Similarly,
Adversary should avoid $\Qmin$ jumping over $\Qmax$ and
thus resetting $S$. They may do so by choosing the initial gap between $Q_i$ and $Q_j$
to be sufficiently large.

Consider the evolution of $\abs{S}$ and of the Boolean value of the statement ``$w[\Qmin]\in k$''
in Adversary's optimal strategy described above. They form a pair $(s,b)$
undergoing a random walk on the domain
$\{0,1,\dotsc,k\}\times \{\stout,\stin\}\setminus\{(0,\stin),(k,\stout)\}$
with the initial state $(s_0,b_0)$, and the following transition rule:
\begin{itemize}
\item If the current state is $(0,\stout)$, then
  \begin{itemize}
  \item \starstep with probability $1/k$, the next state is $(k,\stin)$,
  \item \starstep with probability $1-1/k$, the next state is $(0,\stout)$.
  \end{itemize}
\item If the current state is $(s,\stout)$, for $s\geq 1$, then
  \begin{itemize}
  \item with probability $1/k$, the next state is $(s,\stin)$,
  \item with probability $1/k$, the next state is $(k,\stin)$,
  \item with probability $1-2/k$, the next state is $(s,\stout)$.
  \end{itemize}
\item If the current state is $(s,\stin)$, for $s\geq 1$, then
  \begin{itemize}
  \item \starstep with probability $1/ks$, the next state is $(k,\stin)$,
  \item \starstep with probability $(1/s)(1-1/k)$, the next state is $(0,\stout)$,
  \item with probability $(1-1/s)(1/k)$, the next state is $(s-1,\stin)$,
  \item with probability $(1-1/s)(1/k)$, the next state is $(k,\stin)$,
  \item with probability $(1-1/s)(1-2/k)$, the next state is $(s-1,\stout)$.
  \end{itemize}
\end{itemize}

Here, \starstep indicates that the respective turn is counted by $\overline{B}(s_0,b_0,L)$.

Denote this Markov chain by~$M$.
The chain $M$ is aperiodic because for the state $(0,\stout)$ there is a positive probability of remaining in the state.
The chain is also clearly strongly connected.

Denote by $\pi$ the stationary distribution of $M$. If we start the chain from the stationary distribution, then
$
  \E_{(s,b)\sim \pi} \overline{B}(s,b,L)=L\cdot \Pr_{(s,b)\sim \pi}[\text{\starstep encountered on transition from }(s,b)].
$ 
As chain $M$ is aperiodic, it follows that, for every initial condition $(s,b)$,
\begin{equation}\label{eq:trivialconvergence}
\lim_{L\to\infty} \E [\overline{B}(s,b,L)]/L=\Pr_{(s,b)\sim \pi}[\text{\starstep encountered}].\end{equation}
It remains to compute the probability on the right hand side, and to bound the rate of convergence.

It is routine to verify that the stationary distribution is given by
\[
\Pr_{\pi}[(s,b)]=\begin{cases}
  \frac{1}{k^2 2^{k-1}}\cdot s2^{s-1}&\text{if }b=\stin\text{ and }s\in[k],\\
  \frac{1}{k^2 2^{k-1}}\cdot s2^{s-1}(k-2)&\text{if }b=\stout\text{ and }s\in [k-1],\\
  \frac{1}{k^2 2^{k-1}}\cdot (2^k-1)(k-1)&\text{if }(s,b)=(0,\stout).
\end{cases}
\]
From this we obtain
\[
  p\eqdef \Pr_{(s,b)\sim \pi}[\text{\starstep encountered}]=\frac{1}{k^2 2^{k-1}}\cdot\left((2^k-1)(k-1)+\sum_{s=1}^k 2^{s-1} \right)=\frac{2^k-1}{k 2^{k-1}}.
\]

Let $X_i$ be the characteristic random variable of the event that we visit \starstep at the $i$'th step of the chain.
Then $\overline{B}(s,b,L)=X_1+X_2+\dotsb+X_L$.
To bound the rate of convergence in \eqref{eq:trivialconvergence} we require a tail bound for this sum.
Though there are a number of such results in the literature \cite{lezaud,CLLMchernoff,wagner,dinwoodie,kahale}, the explicit
bounds for non-reversible chains are complicated to state. The simple form of $M$ allows us to give a short and self-contained argument.
The resulting bound is significantly weaker, but it is sufficient for our application.\smallskip

We begin by noting that, from every state, the probability of making a transition to the state $(k,\stin)$ is $1/k$.
Suppose we have two copies of the chain with different starting states. We couple their evolution as follows.
At each step we toss a biased coin that lands heads with probability $1/k$. If it lands on heads, we transition
to $(k,\stin)$ simultaneously in both chains. If it lands on tails, we make steps in the two chains independently conditioned on not transitioning to $(k,\stin)$.
Since the probability that the two chains are in different states after $t$ steps is at most $(1-1/k)^t$, 
this implies, via \cite[Theorem 5.2]{levin_peres}, that the total variation distance to $\pi$ from $M$ after $t$ steps
is at most $(1-1/k)^t$. Let $T$ be a natural number to be chosen later, and partition the sequence
$X_1,\dotsc,X_L$ into $T$ subsequences, each of which is made of samples $T$ steps apart. Let $X_i,X_{T+i},\dotsc,X_{mT+i}$ be any such subsequence,
where $m=m(i)=\lfloor (L-i)/T\rfloor$ and $i\in \{1,2,\dotsc,T\}$. Drop the first element $X_i$; the total variation distance between $(X_{T+i},\dotsc,X_{mT+i})$ and
$m$ independent samples from $\pi$ is at most $m(1-1/k)^T$. Since $X_i$'s are $\{0,1\}$-valued, it follows from the usual Chernoff bound \cite[Theorem A.1.4]{alon_spencer} that
\[
  \Pr[X_{T+i}+\dotsb+X_{mT+i}-mp\geq \lambda]\leq e^{-2\lambda^2/m}+m(1-1/k)^T\leq e^{-2T\lambda^2/L}+me^{-T/k}.
\]
Setting $\lambda=\sqrt{L/2k}$, and using the union bound over $i$ we obtain
\begin{align*}
  \Pr[X_{T+1}+&X_{T+2}+\dotsb+X_L-(L-T)p\geq T\sqrt{L/2k}]\\&\leq \sum_{i\in [T]} \Pr[X_{T+i}+\dotsb+X_{m(i)T+i}-m(i)p\geq \sqrt{L/2k}]\\&\leq (T+L)e^{-T/k}.
\end{align*}
We choose $T=(2kL)^{1/2}p$. Since $X_1+\dotsb+X_T\leq T$ and $T\leq pL$ (as $L>k^3$), we deduce that
\[
  \Pr[\overline{B}(s,b,L)\geq 3pL]\leq 2Le^{-p(2L/k)^{1/2}}.
\]
Since $1/k\leq p\leq 2/k$, the desired bound \eqref{eq:trivialpairtail} follows from an application of \eqref{eq:trivialstrategyarg} 
to the characteristic function of the event ``$\overline{B}(s_0,b_0,T)\geq 3pL$''.

\section{Asymptotics for \texorpdfstring{$\E[P_d-P_0]$}{E[P\textunderscore d-P\textunderscore 0]}}\label{sec:asympt}
In this section we prove \Cref{thm:driftasympt}, which is the asymptotic
\begin{equation}\label{eq:driftasympt}
  \E[P_{\hl}(L)-P_0(L)]=2\sqrt{\hl L/k} \cdot \Bigl(1+O\Bigl(\frac{1}{\hl^{2/3}}+\frac{\log L}{(L/\hl k)^{1/2}}+\frac{\hl^{3/2}}{k^{1/2}}+\frac{k^{1/2}\hl^{3/2}L^{3/2}}{\exp(L^{1/2}k^{-3/2})}\Bigr)\Bigr).
\end{equation}
By \Cref{lem:crudedriftupper} we see that 
\[
\E[P_{\hl}(L)-P_0(L)] \le \hl\sqrt{2L/k}+\hl = 2\sqrt{\hl L/k} \cdot O\bigl(\sqrt{\hl}+\sqrt{\hl k/L}\bigr). 
\]
Let $k_0$ be a constant to be determined.
Note that $\sqrt{\hl} = O\bigl(\frac{k^{1/2}\hl^{3/2}L^{3/2}}{\exp(L^{1/2}k^{-3/2})} + \frac{\hl^{3/2}}{k^{1/2}}\bigr)$ if $L \le k^3$ or $\hl \ge \tfrac{1}{100}k^{1/2}$ or $k\leq C$. Since $\sqrt{\hl k/L} = O\bigl(\frac{\log L}{(L/\hl k)^{1/2}}\bigr)$, in the proof of \eqref{eq:driftasympt} we may assume that $L\geq k^3$, and $\hl\leq \tfrac{1}{100} k^{1/2}$, and $k\geq k_0$.

Recall that, for a word $w$, $\LNDS(w)$ is the length of the longest nondecreasing subsequence in $w$. The study of $\LNDS(w)$ for a
random word began with the case when $w$ is a random permutation. In this case, the asymptotics
was obtained in \cite{logan_shepp} and \cite{versik_kerov}, and then much refined later in 
\cite{baik_deift_johansson_lis}. Similar results for $\LNDS(w)$ itself was first obtained by \cite{tracy_widom}, 
building on the earlier work \cite{baik_deift_johansson_lis} for random permutations. The 
results have been extended and refined in a number of subsequent works, including \cite{johansson,kuperberg}.

The main result of \cite{tracy_widom} shows that, after suitable normalization, $\LNDS(w)$ converges
in distribution to the Tracy--Widom distribution $F_2$. For our purposes, convergence in distribution is insufficient,
as we need an estimate on $\E[\LNDS(w)]$. In \Cref{sec:elnds} we combine the existing results into the following
estimate.
\begin{lemma}\leavevmode\label{lem:elnds}
  \begin{enumerate}[label=(\roman*), ref=(\roman*)]
  \item \label[part]{elnds:constantlength}
For a uniform random word $w\sim [k]^n$,
\[
  \frac{\E[\LNDS(w)]-n/k}{2\sqrt{n}}=1+O\Bigl(\frac{1}{k^{2/3}}+\frac{k^2+k\log n}{n^{1/2}}\Bigr).\hskip12.00093em 
\]
\item \label[part]{elnds:binomiallength}
There exists an absolute constant $C>0$ such that the following holds.
Pick $m\sim \Binom(n,p)$, and then choose $w\sim [k]^m$ uniformly. Then
\[
  \frac{\E[\LNDS(w)]-pn/k}{2\sqrt{pn}}=1+O\Bigl(\frac{1}{k^{2/3}}+\frac{k^2+k\log pn}{(pn)^{1/2}}\Bigr)\qquad \text{whenever }p\frac{n}{\log n}\geq C.
\]
\end{enumerate}
\end{lemma}

We shall analyze the concatenation $A[0]A[1]\dotsb A[L-1]$ appearing in \Cref{prop:lcstolnds} by partitioning the word $A[0]A[1]\dotsb A[L-1]$ into two parts: those subwords $A[\ell]$ that come from trivial
expectant partitions, and those that come from non-trivial expectant partitions.
To analyze the first part we will appeal to
the known results on LNDS in random words (encapsulated in \Cref{lem:elnds}). We will then argue that even
if the remaining partitions $\mathcal{A}[\ell]$ are chosen adversarily, the $\LNDS$ is unlikely to change much.
The adversarial argument is captured by the following Markov game.

\begin{game}\label{game:shuffle} 
 \gamestate Sequence of sets $A[0],\dotsc,A[\ell-1]\subseteq \{0,1,\dotsc,\hl\}$, each of which is
        marked either as `tampered' or `untampered'. 
 \gamestart The sequence is empty.
 \gameduration{$L$}
 \gameturn 
 \begin{itemize}[leftmargin=*]
 \item First, Adversary decides whether to intervene this turn. If they decide to intervene, then they choose a partition $\mathcal{A}$ of $\{0,1,\dotsc,d\}$ into $k$ parts, some of which are possibly empty, and also choose a position where to insert the next set into the list $A[0],\dotsc,A[\ell-1]$ that must be after the last tampered set (if any).
 If they decide not to intervene, then $\mathcal{A}$ becomes the trivial partition.
 \item Then Fortune picks a set $A\in \mathcal{A}$ uniformly among the $k$ sets in $\mathcal{A}$. 
 The set $A$ is then inserted at the chosen position (if Adversary intervened) or at the end of the list (if Adversary abstained this turn). 
 If $\mathcal{A}$ was chosen by Adversary, the inserted set is marked as `tampered'; otherwise, it is `untampered'.
 \end{itemize}
 \gameobjective Adversary aims either to maximize or to minimize $\E[\LNDS(A[0]A[1]\dotsc A[L-1])]$.
 \gamepenalty Adversary pays penalty of $L$ to the objective function if the total number of interventions exceeded $6\hl^2L/k$.
\end{game}

It is not hard to extract from \Cref{prop:update} a Markov chain admitting a morphism to \Cref{game:shuffle}. 
Indeed, we may think of the state of Markov chain from \Cref{prop:update} at time $\ell$ as consisting of
the word $w$, the prefix $w'_{<\ell}$ as well as
all partitions $\mathcal{A}[0],\dotsc,\mathcal{A}[\ell-1]$ and
choices $A[0],\dotsc,A[\ell-1]$. We can then break each step into parts: 
computing $\mathcal{A}[\ell]$ in the first step, and examining $w'[\ell]$ and choosing $A[\ell]$ 
in the second step. By giving Adversary the power to choose the partitions in the even-numbered steps,
and also the power to insert the partitions not only at the end, but also in the middle of the list, we obtain the game above.

\paragraph{Analysis of \texorpdfstring{\Cref{game:shuffle}}{Game \ref{game:shuffle}}.} Since $\LNDS(A[0]\dotsb A[L-1])$ is always between $0$ and $L$, in the optimal play Adversary must never use more than
$6d^2L/k$ interventions. Indeed, any strategy that sometimes uses more than $6d^2L/k$ interventions can be improved
by replacing intervention steps that exceed the $6d^2L/k$ threshold by any other steps (similarly
to the argument in the analysis of \Cref{game:twosimple}).  In addition, since Adversary may insert partitions anywhere in the list, they may as well hold off the interventions until
the very end. Finally, since Adversary is allowed to select trivial partitions, they may as well use all $6d^2L/k$ interventions.

So, under the optimal strategy Adversary first waits for Fortune to
generate a random word for the first $L-6\hl^2L/k$ turns, and then intervenes 
for each of the last $6\hl^2L/k$ turns. Note that because $\hl\leq \tfrac{1}{100} k^{1/2}$,
it follows that $6\hl^2L/k<L$.

Let $\wFort$ be the word generated by Fortune in the first $L'\eqdef L-6\hl^2L/k$ turns. Note that
the length of $\wFort$ might be strictly less than $L'$ because the trivial partitions contain only $\hl+1$ non-empty parts among $k$ parts.
So, $\wFort$ is a uniform random word over the alphabet $\{0,1,\dotsc,\hl\}$ of length $\Binom(L',(\hl+1)/k)$. 
\Cref{lem:elnds}\ref{elnds:binomiallength} tells us that 
\[
  \frac{\E[\LNDS(\wFort)]-L'/k}{2\sqrt{L'(\hl+1)/k}}=1+O\Bigl(\frac{1}{\hl^{2/3}}+\frac{\hl^2+\hl\log pL'}{(pL')^{1/2}}\Bigr)\qquad \text{ whenever }\frac{\hl}{k}\cdot\frac{L'}{\log L'}\geq C,
\]
where $p=(\hl+1)/k$.
We may drop the condition $\frac{\hl}{k}\cdot\frac{L'}{\log L'}\geq C$, as it is satisfied when $L\geq k^3$ and $\hl\leq \tfrac{1}{100}k^{1/2}$ and
$k\geq k_0$ if we define $k_0\eqdef C+10$ with the constant $C$ given by \Cref{lem:elnds}\ref{elnds:binomiallength}.

By changing the variable from $L'$ to $L$, we may rewrite this more conveniently as 
\begin{equation}\label{eq:lndsfortune}
  \frac{\E[\LNDS(\wFort)]-L'/k}{2\sqrt{\hl L/k}}=1+O\Bigl(\frac{1}{\hl^{2/3}}+\frac{\log \hl L/k}{(L/\hl k)^{1/2}}+\frac{\hl^2}{k}\Bigr).
\end{equation}

The following lemma implies that Adversary's intervention cannot decrease $\E[\LNDS_{\hl}-\LNDS_0]$.
\begin{lemma}\label{lem:lndspartitionupper}
Let $w,w'$ be any two words over alphabet $\{0,1,\dotsc,\hl\}$. Then for every partition $\mathcal{A}$ of $\{0,1,\dotsc,\hl\}$ into $k$ (not necessarily nonempty) parts,
\[
  \E_{A\sim \mathcal{A}} \bigl[\LNDS(wAw')\bigr]\geq \LNDS(ww')+1/k.
\]
\end{lemma}
\begin{proof}
Let $uu'$ be a non-decreasing subsequence in $ww'$ of length $\LNDS_{\hl}(ww')$, where $u$ and $u'$ are subsequences in $w$ and $w'$ respectively.
Suppose first that $u$ is non-empty. Let $s\in \{0,1,\dotsc,d\}$ be the last symbol in $u$, and note that
if $s\in A$, then $wAw'$ contains subsequence $usu'$, which is of length $\LNDS_{\hl}(ww')+1$. Hence, $\E\bigl[\LNDS_{\hl}(wAw')\bigr]\geq \LNDS_{\hl}(ww')+1/k$ in this case.
If $u$ is empty, but $u'$ is non-empty, then we let $s$ to be the first symbol of $u'$, and use the same argument to reach the same conclusion.
The remaining case, when both $u$ and $u'$ are empty, is trivial as well. 
\end{proof}

Let $\Penalty$ be the penalty function from \Cref{game:shuffle}. From \Cref{prop:lcstolnds}
and the adversarial game argument it follows that
\[
  \E[P_{\hl}(L)+\Penalty]\geq \inf_{S\in\DStrat(\text{\Cref{game:shuffle}})} \E[\LNDS(A[0]A[1]\dotsb A[L-1])+\Penalty].
\]
Hence, \eqref{eq:lndsfortune} and \Cref{lem:lndspartitionupper} together imply that
\begin{equation}\label{eq:pdlower}
\begin{aligned}
  \E[P_{\hl}(L)+\Penalty]&\geq \E[\LNDS(\wFort)]+\frac{6\hl^2L/k}{k}\\
&= L/k+2\sqrt{\hl L/k} \cdot \Bigl(1+O\Bigl(\frac{1}{\hl^{2/3}}+\frac{\log L}{(L/\hl k)^{1/2}}+\frac{\hl^2}{k}\Bigr)\Bigr).
\end{aligned}
\end{equation} 

We next tackle the upper bound on $\E[P_d(L)]$. We will use the following simple fact.
\begin{lemma}\label{lem:lndspartitionsub}
Suppose a word $w$ (which we view as a sequence of symbols) is partitioned
into two subsequences $u_1$ and $u_2$. 
Then
$
  \LNDS(w)\leq \LNDS(u_1)+\LNDS(u_2).
$
\end{lemma}
\begin{proof}
The restrictions of a nondecreasing subsequence in $w$ to $I_1$ and $I_2$ are nondecreasing subsequences
in $u_1$ and $u_2$, respectively. So, $\LNDS(w)\leq \LNDS(u_1)+\LNDS(u_2)$ follows.
\end{proof}

Because in \Cref{game:shuffle} Adversary can select the insertion position only
after the already-tampered sets, this reduces our task to analyzing the following Markov game:
\begin{game}\label{game:lnds}
\gamestate Finite word $w$ over alphabet $\{0,1,\dotsc,\hl\}$.
\gamestart The word $w$ is empty.
\gameduration{$L$}
\gameturn Adversary chooses a partition $\mathcal{A}$ of $\{0,1,\dotsc,\hl\}$ into $k$ parts, some of which are possibly empty.
Then Fortune picks a set $A\in\mathcal{A}$ uniformly among the $k$ sets in $\mathcal{A}$. The elements of $A$ are then
appended in descending order to~$w$.
\gameobjective Adversary aims to maximize $\E[\LNDS(w)]$.
\end{game}
Thanks to the connection between $\LNDS$ and the evolution of $P_0,P_1,\dotsc,P_{\hl}$ from \Cref{prop:lcstolnds},
we see that this game is in fact equivalent to \Cref{game:particles}, for which we already gave an upper bound
in \Cref{lem:adversariallndsupper}. Combining that bound with the upper bound from \Cref{lem:lndspartitionsub}
we obtain
\begin{align*}
  \E[P_{\hl}(L)-\Penalty]&\leq \sup_{S\in\DStrat(\text{\Cref{game:shuffle}})} \E[\LNDS(A[0]A[1]\dotsb A[L-1])-\Penalty]+\hl\\
                        &\leq \E[\LNDS(\wFort)]+
\sup_{\substack{S\in\DStrat(\text{\Cref{game:lnds}})\\\text{for }6\hl^2L/k\text{ turns}} }\E[\LNDS(A[0]\dotsb A[6\hl^2L/k-1])]+\hl\\
                        &\leq L'/k+2\sqrt{\hl L/k}\left(1+O\Bigl(\frac{1}{\hl^{2/3}}+\frac{\log \hl L/k}{(L/\hl k)^{1/2}}+\frac{\hl^2}{k}\Bigr)\right)\\
                        &\qquad       + \frac{6d^2L}{k^2} + \hl\sqrt{2(6d^2L/k)/k}+2\hl\\
                        &\leq \frac{L}{k}+2\sqrt{dL/k} \cdot \Bigl(1+O\Bigl(\frac{1}{\hl^{2/3}}+\frac{\log dL/k}{(L/\hl k)^{1/2}}+\frac{\hl^{3/2}}{k^{1/2}}\Bigr)\Bigr).
\end{align*}

\Cref{lem:trivialpartitionbound} implies that $\E[\Penalty]\leq 4\hl^2L^2\exp(-L^{1/2}k^{-3/2})$. So,
from \eqref{eq:pdlower} and the preceding bound we deduce that
\[
  \E[P_{\hl}(L)]=\frac{L}{k}+2\sqrt{\hl L/k} \cdot \Bigl(1+O\Bigl(\frac{1}{\hl^{2/3}}+\frac{\log dL/k}{(L/\hl k)^{1/2}}+\frac{\hl^{3/2}}{k^{1/2}}+\frac{\hl^2L^{3/2}}{\exp(L^{1/2}k^{-3/2})}\Bigr)\Bigr).
\]
Hence \Cref{thm:driftasympt} follows since $\E[P_0(L)]=L/k$.

\section{Remarks and open problems}
\begin{itemize}
\item  We believe that \Cref{lem:adversariallndsupper} should hold with $C\sqrt{dL/k}+d$ instead of $d\sqrt{2L/k}+d$.
Any upper bound in \Cref{lem:adversariallndsupper} that is sublinear in $d$ would imply that $\lim \gamma_k'=\tfrac{1}{4}$
(assuming that constants $\gamma_k'$ exist). That is because we can upper bound $\E[P_d-P_0]$ by smallest
of the bound in \Cref{lem:crudedriftupper} and the bound in such improved \Cref{lem:adversariallndsupper}.
That can then be used in an argument similar to that in \Cref{sec:upperbound}.

\item The formalism of Markov games can be generalized. Most notably instead of giving an adversary, for each state $x_n$,
a choice of states $R_n(x_n)\subset \Omega_{n+1}$ where they can take the chain to, we give them a choice of a probability distributions
on $\Omega_{n+1}$. Perhaps the most natural case is when the set of allowable probability distributions is a convex set
in the space of probability distributions on $\Omega_{n+1}$. Doing so would eliminate the need for the set $I$ in the definition
of a Markov game. Since we did not need the extra generality, we opted for less abstract presentation.

\end{itemize}

\bibliographystyle{plain}
\bibliography{unbalancedlcs}

\appendix
\section{Proofs of \texorpdfstring{\Cref{prop:update,prop:lcstolnds}}{Propositions \ref{prop:update} and \ref{prop:lcstolnds}}}\label{sec:props}
In this appendix, we supply proofs that were omitted in the introduction.

\paragraph{Proof of \texorpdfstring{\Cref{prop:update}}{Proposition \ref{prop:update}}.}
\textit{Case $w[P_i(L)]=w'[L]$:} Let $u$ be a common subsequence of $w_{<P_i(L)}$ and $w'_{< L}$
that is obtained from $w_{< P_i(L+1)}$ by removing at most $i$ symbols.
By appending $w'[L]$ to $u$ we obtain a subsequence of $w'_{< L+1}$
that is a witness to the inequality $P_i(L+1)\geq P_i(L)+1$. As the reverse inequality
$P_i(L+1)\leq P_i(L)+1$ always holds, this shows that $P_i(L+1)=P_i(L)+1$.

\textit{Case $w[P_i(L)]\neq w'[L]$:} 
Note that we always have $P_i(L+1)\geq \max\bigl(P_i(L),P_{i-1}(L+1)+1\bigr)$. We must show the reverse inequality.

Suppose $P_i(L+1)>P_i(L)$, and so $P_i(L+1)=P_i(L)+1$. Let $\tilde{u}$ be the common subsequence of $w'_{< L+1}$ and $w_{<P_i(L+1)}$ obtained from the latter by
omitting $i$ symbols.  Since $w[P_i(L)]\neq w'[L]$, the last symbol in $\tilde{u}$ is either different from $w[P_i(L)]$ or from $w'[L]$.
In the former case, $P_i(L+1)\leq P_{i-1}(L+1)+1$ holds, whereas in the latter case $P_i(L+1)\leq P_i(L)$ holds.

\paragraph{Proof of \texorpdfstring{\Cref{prop:lcstolnds}}{Proposition \ref{prop:lcstolnds}}.} The proof is by induction on $L$ and~$i$. When $L=0$, we have $P_i(0)=i$ for all $i$
and the only non-decreasing subsequence of the empty word is empty. Similarly, if $i=0$, the claim is trivially true. 
Suppose now that we wish to prove the claim for the pair ($L+1$,$i$), and that the claim holds for 
all smaller pairs. 

Suppose $i\in A[L]$. In this case, $P_i(L+1)=P_i(L)+1$ according to \Cref{prop:update}. On the other hand,
we may append $i$ to any non-decreasing subsequence in $A[0]\dotsb A[L-1]$ using only symbols from $\{0,1,\dotsc,i\}$,
and so $\LNDS_i(A[0]\dotsb A[L-1]A[L])=\LNDS_i(A[0]\dotsb A[L-1]A[L])+1$.

Suppose $i\notin A[L]$. In this case, consider the longest non-decreasing sequence in $A[0]\dotsb A[L]$
that uses only symbols from $\{0,1,\dotsc,i\}$. If the sequence contains no $i$, then 
$\LNDS_i(A[0]\dotsb A[L])=\LNDS_{i-1}(A[0]\dotsb A[L])=(P_{i-1}(L+1)-i)+1$. If the sequence
contains $i$, then it does not contain any symbols from $A[L]$, in which case
$\LNDS_i(A[0]\dotsb A[L])=\LNDS_{i-1}(A[0]\dotsb A[L-1])=P_i(L)-i$. Both cases match the behavior of $P_i(L+1)$ from \Cref{prop:update}.

\section{Expected modulus of a simple random walk}\label{sec:randomwalk}
Given $T \in \N$. Consider a particle doing a simple random walk starting from $1/2$ for $T$ steps. Let $\Delta'$ be the final position of this particle. 

\begin{proposition}\label{prop:randomwalk}
	We have that 
	\[
	\E\bigl[\abs{\Delta'}\bigr]=
	\begin{dcases}
		2^{-T} \cdot \Bigl(T+\frac{1}{2}\Bigr)\binom{T}{T/2},&\text{if }T\text{ is even},\\
		2^{-T}\cdot 2T\binom{T-1}{(T-1)/2},&\text{if }T\text{ is odd}.
	\end{dcases}
	\]
\end{proposition}

\begin{proof}
        Suppose $T$ is even, and write it as $T = 2m$. Set $p \eqdef 2^{-2m}$. 
	
	The binomial theorem tells us that $\Pr[\Delta' = 1/2] = p \cdot \binom{2m}{m}$, and that 
	\[
	\Pr[\Delta' = 1/2 + 2r] = \Pr[\Delta' = 1/2 - 2r] = p \cdot \tbinom{2m}{m+r} = p \cdot \tbinom{2m}{m-r}
	\]
	for $r = 1, \dotsc, m$. Then
	\begin{align*}
	p^{-1} \E\bigl[|\Delta'|\bigr] &= \tfrac{1}{2} \cdot \tbinom{2m}{m} + \sum_{r=1}^{m} 2r \cdot 2 \cdot \tbinom{2m}{m+r} \\
	&= \tfrac{1}{2} \cdot \tbinom{2m}{m} + 4 \cdot \sum_{r \ge 1} (m+r) \cdot \tbinom{2m}{m+r} - 4m \cdot \sum_{r \ge 1} \tbinom{2m}{m+r} \\
	&= \tfrac{1}{2} \cdot \tbinom{2m}{m} + 4 \cdot \sum_{r \ge 1} 2m \cdot \tbinom{2m-1}{m+r-1} - 4m \cdot \sum_{r \ge 1} \tbinom{2m}{m+r} \\
	&= \tfrac{1}{2} \cdot \tbinom{2m}{m} + 8m \cdot \tfrac{1}{2} \cdot 2^{2m-1} - 4m \cdot \tfrac{1}{2} \cdot \bigl(2^{2m}-\tbinom{2m}{m}\bigr) \\
	&= (2m+\tfrac{1}{2}) \cdot \tbinom{2m}{m}. 
	\end{align*}
Hence $\E\bigl[|\Delta'|\bigr] = 2^{-T}(T+\frac{1}{2})\binom{T}{T/2}$. 

Suppose $T$ is odd, and write it as $T = 2m+1$. Set $p \eqdef 2^{-(2m+1)}$. 

The binomial theorem tells us that 
\[
\Pr[\Delta' = 1/2 + 2r - 1] = \Pr[\Delta' = 1/2 - 2r + 1] = p \cdot \tbinom{2m+1}{m+r} = p \cdot \tbinom{2m+1}{m+1-r}
\]
for $r = 1, \dotsc, m+1$. Then
\begin{align*}
	p^{-1} \E\bigl[|\Delta'|\bigr] &= \sum_{r=1}^{m} (2r-1) \cdot 2 \cdot \tbinom{2m+1}{m+r} \\
	&= 4 \cdot \sum_{r \ge 1} (m+r)\tbinom{2m+1}{m+r} - (4m+2) \cdot \sum_{r \ge 1} \tbinom{2m+1}{m+r} \\
	&= 4 \cdot \sum_{r \ge 1} (2m+1) \cdot \tbinom{2m-1}{m+r-1} - (4m+2) \cdot \sum_{r \ge 1} \tbinom{2m+1}{m+r} \\
	&= (8m+4) \cdot \tfrac{1}{2} \cdot \bigl(2^{2m}+\tbinom{2m}{m}\bigr) - (4m+2) \cdot \tfrac{1}{2} \cdot 2^{2m+1} \\
	&= (4m+2) \tbinom{2m}{m}. 
\end{align*}
Hence $\E\bigl[|\Delta'|\bigr] = 2^{-T+1}T\binom{T-1}{(T-1)/2}$.
\end{proof}

\section{Behavior of LNDS for growing alphabet}\label{sec:elnds}
Here, we prove \Cref{lem:elnds} describing the behavior of $\E[\LNDS(w)]$ for a random $w\sim [k]^n$ uniformly in $k$ and~$n$.
A closely related work is \cite{breton_houdre}, which asserts, as a special case, convergence of $\LNDS(w)$ to the Tracy--Widom
distribution as both $k\to\infty$ and $n\to\infty$, subject to appropriate growth conditions on $k$. Sadly, we 
are unable to use \cite{breton_houdre} because it does not claim convergence of moments. Also, the proof in \cite{breton_houdre}
has a minor gap --- after the application of the one-dimensional strong approximation result of Sakhanenko
to several dependent random variables, the resulting Brownian motions will be dependent, but nothing
is proved about their dependency. In a private correspondence, the authors of \cite{breton_houdre} indicated 
that this gap can be fixed by appealing to results in \cite{houdre_litherland} and \cite{houdre_xu}. We take 
this opportunity to present an alternative argument. An advantage of our argument is that we obtain an explicit error 
term.

\paragraph{Poissionization and de-Poissonization.} 
Denote the Poisson random variable with mean $\lambda$ by $\Pois(\lambda)$, and define random variables
\begin{align*}
  L_n&\eqdef \LNDS(w)-n/k\text{ for }w\sim [k]^n,\\
  \widetilde{L}_{\lambda}&\eqdef L_{\Pois(\lambda)}.
\end{align*}
The variable $\widetilde{L}_{\lambda}$ serves a standard purpose: it is substantially easier to analyze than $L_n$.
We shall recover $\E L_n$ from $\E \widetilde{L}_{\lambda}$ via de-Poissonization argument in the following lemma. 
\begin{lemma} \label{lem:depoissonization}
Suppose $s_0,s_1,\dotsc$ is an increasing sequence satisfying $s_m-s_{m-1}\leq m^A$ with $A\geq 0$.
Define $\widetilde{s}_{\lambda} \eqdef s_{\Pois(\lambda)}$.
\newcommand*{\WidestLHS}{\abs{\E[s_{\Binom(n,p)}]-\widetilde{s}_{\lambda}}}
\newcommand*{\WidestCS}{p\E[\widetilde{s}_{\lambda_+}-\widetilde{s}_{\lambda_-}]+13/n^A,}
\newcommand*{\WidestRHS}{p\frac{n}{\log n}\geq n_0(A).}
\newcommand*{\FormatLHS}[1]{\makebox[\widthof{$\WidestLHS$}][r]{$#1$}}%
\newcommand*{\FormatCS}[1]{\makebox[\widthof{$\WidestCS$}][l]{$#1$}}%
\newcommand*{\FormatRHS}[1]{\makebox[\widthof{$\WidestRHS$}][l]{$#1$}}%
\begin{enumerate}[label=(\roman*), ref=(\roman*)]
	\item \label{depois:simple} 
	We have 
	\[
	\FormatLHS{\E\widetilde{s}_{\lambda^-}-4/n^A} \leq \FormatCS{s_n\leq \E\widetilde{s}_{\lambda^+}+4/n^{A},}\qquad\text{whenever }\FormatRHS{n\geq n_0(A),}
	\]
	where $\lambda^{\pm}=n\pm(A+1)\sqrt{2n \log n}$.
	\item \label{depois:binomial}
	We also have
	\[
	\FormatLHS{\abs{\E[s_{\Binom(n,p)}-\widetilde{s}_{pn}]}}\leq \FormatCS{p\E[\widetilde{s}_{\lambda_+}-\widetilde{s}_{\lambda_-}]+13/n^A,}\qquad \text{whenever }\FormatRHS{p\frac{n}{\log n}\geq n_0(A).}
	\]
	where $\lambda^{\pm}=pn\pm(2A+3)\sqrt{2pn \log n}$.
\end{enumerate}
\end{lemma}
An often-quoted result that is similar to \ref{depois:simple} appears in \cite[Lemma 2.5]{johansson_unitary}; 
a slightly more general and streamlined version is reproduced in \cite[Lemma~2.31]{romik}. 
The part \ref{depois:binomial}, as far as we are aware, is new. Without the factor $p$, it follows
easily from \ref{depois:simple} and the Chernoff bounds. The extra factor allows for simpler application, 
and results in a slightly stronger error term in~\Cref{thm:driftasympt}.

\begin{proof} We shall use the tail bound for the Poisson distribution from \cite{cannone_poisson}:
\begin{equation}\label{eq:poissonconc}
  \Pr\bigl[\abs{\Pois(\lambda)-\lambda}\geq x\bigr]\leq 2e^{-\frac{x^2}{2(\lambda+x)}},\qquad x\geq 0.
\end{equation}
\paragraph{Proof of \protect\ref{depois:simple}.}
Let $\lambda=\lambda^-$. Then
\begin{equation}\label{eq:poissonuppertail}
  \E\widetilde{s}_{\lambda}\leq s_n+\sum_{m>n} (s_m-s_{m-1})\Pr[\Pois(\lambda)\geq m]\leq s_n+2\sum_{m>n} m^A e^{-\frac{(m-\lambda)^2}{2m}}.
\end{equation}
The ratio between consecutive terms in the sum is
$
  (1+1/m)^A\cdot \exp(-\frac{1}{2}+\frac{\lambda^2}{2m(m+1)})
$, which is decreasing in~$m$. 
Hence, by our choice of $\lambda$ and the condition, the ratio is at most $\exp(-2/n^{1/2})$ for large $n$, and so
$\E\widetilde{s}_{\lambda}\leq s_n+4n^{A+1/2}e^{-(n-\lambda)^2/2n}\leq s_n+4/n^A$.

Similarly, we let $\lambda=\lambda^+$ and compute
\begin{align*}
  \E\widetilde{s}_{\lambda}&\geq s_n-\sum_{m< n} (s_{m+1}-s_m)\Pr[\Pois(\lambda)\leq m]\geq s_n-2\sum_{m<n} m^A e^{-\frac{(\lambda-m)^2}{2\lambda-m}}
\\&\geq s_n-2n^{A+1}e^{-\frac{(\lambda-n)^2}{2\lambda}}\geq s_n-4/n^A.
\end{align*}

\paragraph{Proof of \protect\ref{depois:binomial}.} 
Define $\lambda\eqdef pn$ and consider probability distributions $\Binom(n,p)$ and $\Pois(\lambda)$.
In \mbox{\cite[Theorem~1]{barbour_hall}} it is shown that the total variation distance between the two distributions
satisfies $\dTV\eqdef \dTV(\Binom(n,p),\Pois(\lambda))\leq p(1-\exp(-\lambda))$.
Consider the difference 
\[
  \Delta(m)\eqdef \Pr[\Binom(n,p)=m]-\Pr[\Pois(\lambda)=m],
\]
and note that 
\[
  \E[s_{\Binom(n,p)}]=\E[\widetilde{s}_{\lambda}]+\sum_m s_m \Delta(m)
\]
Let $x\eqdef (A+2)\sqrt{pn\log n}$. When $m\leq \lambda+x$ and $\Delta(m)>0$, we upper bound $s_m$ by $s_{\lambda+x}$. Similarly,
when $m\geq \lambda - x$ and $\Delta(m)<0$, we upper bound $-s_m$ by $-s_{\lambda-x}$. We thus obtain
\begin{align*}
  \sum_m s_m \Delta(m)&\leq \sum_{m:\Delta(m)>0} s_{\lambda+x} \Delta(m)+\sum_{\substack{m:\Delta(m)>0\\m>\lambda+x}} (s_m-s_{\lambda+x})\Delta(m)\\
                      &\qquad +\sum_{m:\Delta(m)<0} s_{\lambda-x}\Delta(m)-\sum_{\substack{m:\Delta(m)<0\\m<\lambda-x}} (s_{\lambda-x}-s_m)\Delta(m)\\
                     &\leq \dTV (s_{\lambda+x}-s_{\lambda-x})+n^{A+1}\Pr[\Binom(n,p)>\lambda+x]+(\lambda-x)^{A+1}\Pr[\Pois(\lambda)<\lambda-x].\\
\intertext{Applying the asymmetric Chernoff bound from \cite[Theorem~A.1.11]{alon_spencer} and \eqref{eq:poissonconc}, we derive}
  \sum_m s_m \Delta(m)&\leq p(s_{\lambda+x}-s_{\lambda-x})+n^{A+1}(e^{-x^2/2\lambda+x^3/2\lambda^2}+2e^{-x^2/2\lambda})\\
                     &\leq p(s_{\lambda+x}-s_{\lambda-x})+3/n^A.
\shortintertext{Similarly, using the asymmetric Chernoff bound from \cite[Theorem~A.1.13]{alon_spencer}, we obtain}
-\sum_m s_m \Delta(m)&\leq \dTV (s_{\lambda+x}-s_{\lambda-x})+n^{A+1}\Pr[\Binom(n,p)<\lambda-x]+\sum_{m>\lambda+x} m^{A+1} \Pr[\Pois(\lambda)=m]\\
                    &\leq \dTV (s_{\lambda+x}-s_{\lambda-x})+1/n^A+\sum_{m>\lambda+x} m^{A+1} e^{-(m-\lambda)^2/2m}.
\end{align*}
This can be bounded as in \eqref{eq:poissonuppertail} to obtain
\[
-\sum_m s_m \Delta(m)\leq \dTV (s_{\lambda+x}-s_{\lambda-x})+5/n^A.
\]
By part \ref{depois:simple}, $s_{\lambda+x}\leq \E\widetilde{s}_{\lambda_+}+4/n^A$ and $s_{\lambda-x}\geq \E \widetilde{s}_{\lambda_-}-4/n^A$.
\end{proof}

To apply the preceding lemma we will need the following fact.
\begin{lemma} \label{lem:lnmonotone}
The expectation $\E[L_n]$ is an increasing function of~$n$.
\end{lemma}
\begin{proof} Think of a random word $w\sim [k]^{n+1}$ as consisting of the first symbol $w[0]$ and the suffix $w'$ of length~$n$.
Let $A\in [k]$ be the largest symbol such that $w'$ contains a non-decreasing subsequence of length $\LNDS(w')$ whose first symbol is $A$.
Then $\E[\LNDS(w)-\LNDS(w')\mid w']=A/k$. Since $A\geq 1$, it follows that $\E[L_{n+1}-L_n\mid w']=(A-1)/k\geq 0$.
\end{proof}

\paragraph{Geometric view.}
We will use a geometric representation for $\widetilde{L}_{\lambda}$: Let $T\eqdef \lambda/k$ and imagine $k$ independent Poisson
point processes $P^{(1)},\dotsc,P^{(k)}$ on the interval $[0,T]$, each is of constant intensity~$1$. For each $i\in [k]$, 
replace points of $P^{(i)}$ with the symbol $i$, and read
the symbols left-to-right to obtain a random word of length $\Pois(\lambda)$.
Then $\widetilde{L}_{\lambda}$ is the normalized length of the $\LNDS$ in this word.

It is convenient to normalize the processes $P^{(1)},\dotsc,P^{(k)}$ by subtracting their means.
We define the processes $Q^{(1)},\dotsc,Q^{(k)}$ by $Q^{(i)}_t\eqdef P^{(i)}_t-t$. Note that
these processes still have independent increments. Then
\[
  \widetilde{L}_{\lambda}=-\frac{1}{k}\sum_{i=1}^k Q^{(i)}_T +\max_{\substack{0=t_0\leq t_1\leq t_2\leq \dotsb\\\dotsb\leq t_{k-1}\leq t_k=T}} \sum_{i=1}^k \bigl(Q^{(i)}_{t_i}-Q^{(i)}_{t_i-1}\bigr).
\]

\paragraph{Strong approximation.} 
We shall approximate the processes $Q^{(1)},\dotsc,Q^{(k)}$ by Brownian motions.
\begin{lemma}\label{lem:strongapprox}
Let $Q$ be as above. Then there are independent standard Brownian motions $B^{(1)},\dotsc,B^{(k)}$
such that
\[
  \Pr\left[\max_{i\in[k]}\max_{t \in [0,T]} |Q^{(i)}_t-B^{(i)}_t|>x\right]\leq k\exp(-cx),\qquad\text{ for }x\geq C\log T,
\]
where $c>0$ is a constant independent of $k$ and~$T$.
\end{lemma}
\begin{proof}
Consider the values of $Q^{(i)}_t$ at integral values of~$t$. 
For a fixed value of $i$, the
increments $Q^{(i)}_1-\nobreak Q^{(i)}_0,\dotsc,Q^{(i)}_T-Q^{(i)}_{T-1}$ are independent copies
of~$\Pois(1)-1$. By the Koml\'os--Major--Tusn\'ady strong approximation theorem \cite{komlos_major_tusnady}
there are constants $c_0,C_0>0$ that depend only on the distribution $\Pois(1)-1$, 
and Brownian motions $B^{(i)}$ such that 
\[
  \Pr[\Delta^{(i)}\geq C_0\log T+x]\leq \exp(-c_0 x),
\]
where $\Delta^{(i)}\eqdef \max_{m=0,1,\dotsc,T} |Q^{(i)}_m-B^{(i)}_m|$.
Note that since $Q^{(1)},\dotsc,Q^{(k)}$ are independent, we may choose $B$'s to be independent as well.

Since the Poisson point process is monotone, for $t\in [m,m+1]$ we have
\[
  Q^{(i)}_t-B^{(i)}_t\leq Q^{(i)}_{m+1}-B^{(i)}_{m+1}+(B^{(i)}_{m+1}-B^{(i)}_t).
\]
Because, as a function of $t$, the difference $B^{(i)}_{m+1}-B^{(i)}_t$ is the standard Brownian
motion, it follows that  
\[
  \Pr[\exists t\in [0,T],\ i\in [k] \text{ s.t. } Q^{(i)}_t-B^{(i)}_t\geq \Delta^{(i)}+x]
       \leq kT \Pr[\max_{t\in[0,1]} B_t \geq x]\leq kT\exp(-x^2/2).
\]
Since $Q^{(i)}_t$ decreases by at most $1$ on any interval of length $1$, we similarly
derive
\[
  \Pr[\exists t\in [0,T],\ i\in [k] \text{ s.t. } Q^{(i)}_t-B^{(i)}_t\leq -\Delta^{(i)}-x-1]
       \leq kT \Pr[\max_{t\in[0,1]} B_t \geq x+1]\leq kT\exp(-x^2/2).  
\]
Putting these together we obtain
\[
  \Pr[\max_{i\in[k]} \max_{t\in [0,T]} \abs{Q^{(i)}_t-B^{(i)}_t}\geq C_0\log T+2x+1]\leq k\exp(-c_0x)+kT\exp(-x^2/2),
\]
from which the promised inequality follows by appropriate choice of constants $c$ and~$C$.
\end{proof}
Let
\[
  R\eqdef -\frac{1}{k}\sum_{i=1}^k B^{(i)}_T +\max_{\substack{0=t_0\leq t_1\leq t_2\leq \dotsb\\\dotsb\leq t_{k-1}\leq t_k=T}} \sum_{i=1}^k \bigl(B^{(i)}_{t_i}-B^{(i)}_{t_i-1}\bigr).
\]
From \Cref{lem:strongapprox} we see that $\Pr[\abs{R-\widetilde{L}_{\lambda}}>(2k+1) x]\leq k \exp(-cx)$ for $x\geq C\log T$.
In particular,
\begin{equation}\label{eq:gaussianerror} 
\begin{aligned}
  \abs[\big]{\E[R-\widetilde{L}_{\lambda}]}&\leq \E\bigl[\abs{R-\widetilde{L}_{\lambda}}\bigr]=(2k+1)\int_0^{\infty} \Pr\bigl[\abs{R-\widetilde{L}_{\lambda}}\geq (2k+1)x\bigr]\,dx\\&\leq 
  (2k+1)\Bigl(C\log T+ \int_{C\log T}^{\infty}k\exp(-cx)\,dx\Bigr)
   =O(k\log T+k^2).
\end{aligned}
\end{equation}

The expectation of the first term in the definition of $R$ vanishes. When it comes to the second term, the scaling invariance of the Brownian motion
implies that it suffices to consider the case $T=1$. In that case, the main result of \cite{baryshnikov} asserts that the maximum in the definition of $R$ is equal in the law to $\lmax$,
the largest eigenvalue of a $k$-by-$k$ Gaussian Unitary Ensemble (GUE). Hence,
\[
  \E[R]=\sqrt{T}\cdot \E[\lmax]=\sqrt{\lambda/k}\cdot \E[\lmax].
\]

We have been unable to find asymptotic for $\E[\lmax]$ in the literature. The most precise result that we are aware of is in \cite{ledoux}: On the bottom of page~27 in \cite{ledoux} it is asserted 
that
\begin{equation}\label{eq:elmax}
  1-\frac{1}{C'k^{2/3}}\leq \frac{\E[\lmax]}{2\sqrt{k}}\leq 1-\frac{1}{Ck^{2/3}}
\end{equation}
for some $C,C'>0$ and all $k\geq 1$ (note that \cite{ledoux} uses a non-standard normalization for GUE, which
we accounted for when copying the result). 
The proof of \eqref{eq:elmax} in \cite{ledoux} relies on the estimate~(2.11) therein, which in turn relies on the tail bound in 
Proposition~2.4, but no proof of the proposition is given\footnote{Proposition~2.4 in \cite{ledoux} is introduced by ``\dots there is no doubt that a similar Riemann--Hilbert analysis might be
performed anagously [sic] for these examples, and that the statements corresponding to Proposition 2.3 hold true. We may for example guess the following\dots''}. 
However, on page~51 of \cite{ledoux} another tail bound is given in (5.16). While it is weaker than Proposition~2.4, it is strong enough to imply~(2.11), and hence also \eqref{eq:elmax} above.

Putting \eqref{eq:gaussianerror} and \eqref{eq:elmax} together yields
\[
  1-\frac{1}{C'k^{2/3}}+O\Bigl(\frac{k^2+k\log \lambda}{\lambda^{1/2}}\Bigr)\leq \frac{\E[\widetilde{L}_{\lambda}]}{2\sqrt{\lambda}}\leq 1-\frac{1}{Ck^{2/3}}+O\Bigl(\frac{k^2+k\log \lambda}{\lambda^{1/2}}\Bigr).
\]
Combining \Cref{lem:lnmonotone} and \Cref{lem:depoissonization}\ref{depois:simple} applied with $A=2$ to $s_n=\E[L_n]$, we obtain
\[
  1-\frac{1}{C'k^{2/3}}+O\Bigl(\frac{k^2+k\log n}{n^{1/2}}\Bigr)\leq \frac{\E[L_n]}{2\sqrt{n}}\leq 1-\frac{1}{Ck^{2/3}}+O\Bigl(\frac{k^2+k\log n}{n^{1/2}}\Bigr),
\]
which proves \cref{elnds:constantlength} of \Cref{lem:elnds}. The \cref{elnds:binomiallength} follows similarly with the help of
\Cref{lem:depoissonization}\ref{depois:binomial}.

\end{document}